\documentclass[12pt,a4paper]{article}

\usepackage[dvips]{graphics}

\usepackage[dvips]{epsfig}
\usepackage{authblk}

\newtheorem{lemma}{Lemma}

\newtheorem{cor}[lemma]{Corollary}
\newtheorem{prop}[lemma]{Proposition}
\newtheorem{defn}{Definition}
\newtheorem{rem}{Remark}

\newcommand{\dimo}{\noindent \emph{Proof. }}
\newcommand{\qed}{\\ \rightline{$\Box$}\\}

\usepackage{latexsym}
\usepackage{amsfonts}
\usepackage{amssymb}
\usepackage{amsmath}
\usepackage[english]{babel}

\def\G{\Gamma}

\begin{document}

\title{COMPUTING MATVEEV'S COMPLEXITY \\ VIA CRYSTALLIZATION THEORY: \\ THE BOUNDARY CASE}

\bigskip

 \author{Maria Rita CASALI - Paola CRISTOFORI}

\affil{Dipartimento di Matematica \\
Universit\`a di Modena e Reggio Emilia \\ Via Campi 213 B, I-41125 MODENA (Italy)}

\maketitle

\begin{abstract}
The notion of  {\it Gem-Matveev complexity} ({\it GM-complexity}) has been introduced within crystallization theory, as a combinatorial method to estimate Matveev's complexity of closed 3-manifolds; it yielded upper bounds for interesting classes of such manifolds.
In this paper we extend the definition to the case of non-empty boundary and prove that for each compact irreducible and boundary-irreducible 3-manifold it coincides with the {\it modified Heegaard complexity} introduced in \cite{[CMV]}.
Moreover, via GM-complexity, we obtain an estimation of Matveev's complexity for all Seifert 3-manifolds with base $\mathbb D^2$ and two exceptional fibers and, therefore, for all torus knot complements.
\end{abstract}

\bigskip

\small{

\thanks{

{\it 2010 Mathemathics Subject Classification:} Primary 57M27, 57N10. Secondary 57M15, 57M25.

\smallskip

 }

\bigskip
\bigskip

\section{\hskip -0.7cm . Introduction}

\bigskip

The {\it Matveev's complexity}  is a well-known invariant for 3-manifolds, defined in \cite{[M$_1$]} as the minimum number of true vertices among all almost simple spines of the manifold.
The 3-sphere, the real projective space, the lens space
$L(3,1)$ and the spherical bundles $\mathbb S^1\times \mathbb S^2$
and $\mathbb S^1\tilde\times \mathbb S^2$ have complexity zero by
definition. Apart from these special cases, the Matveev's complexity $c(M)$ of a closed prime
3-manifold $M$ turns out to be the minimum number of tetrahedra needed to obtain $M$ via face paring of them
(see \cite[Proposition 2]{[M$_1$]}, together with the related Remark).

The idea of estimating $c(M)$ by using Heegaard
decompositions is already suggested in the foundational paper
\cite{[M$_1$]}: if  $H = (S, \mathbf {v}, \mathbf {w})$ is a Heegaard diagram
of $M$, an upper bound for $c(M)$ is easily provided by the almost simple spines of $M$
obtained by adding to the surface $S$ the meridian disks
corresponding to the systems of curves $\mathbf {v}$ and $\mathbf {w}$ and by removing the 2-disk
corresponding to an arbitrary region $\bar R$ of $S \setminus (\mathbf {v} \cup \mathbf {w}).$
Obviously, the true vertices of the considered spine of $M$ are the intersection points of the curves of the two
systems, with the exception of those which lie on the
boundary of the region $\bar R$; hence
$$c(M)\leq n-m,$$
where $n$ (resp. $m$) denotes the number of intersection points between  $\mathbf {v}$ and $\mathbf {w}$ (resp. the number of intersection points contained in $\bar R$).

Starting from this idea, two different approaches to Matveev's
complexity computation have been recently developed. The first
one has been introduced in 2004 for closed 3-manifolds and is based on {\it crystallization theory}\footnote{Basic elements of this theory are recalled in section 3, where a wide bibliography is also indicated.}: it has led to the notion of {\it
Gem-Matveev complexity}, $GM$-complexity for short (see
\cite{[C$_4$]}, together with subsequent papers \cite {[C$_5$]}
and \cite{[CC$_1$]}). Later, in 2010, the {\it modified Heegaard complexity}
($HM$-complexity) of a compact  3-manifold has been
defined via {\it generalized Heegaard diagrams} (see \cite{[CMV]}).\footnote{Actually, in the original paper  \cite{[CMV]} the definition is given only in the orientable case, while the non-orientable extension is due to \cite{[CCM]}.} Both invariants have been proved to be upper bounds for the Matveev's
complexity.

Moreover, the coincidence of Gem-Matveev complexity and
modified Heegaard complexity is established in \cite{[CCM]} for each closed 3-manifold.

The aim of the present paper is to extend the definition of Gem-Matveev complexity to the case of 3-manifolds with non-empty boundary (Section 4), and to prove that $GM$-complexity and $HM$-complexity  turn out to be useful different tools to compute the same upper bound for Matveev's complexity,  in the whole setting of compact irreducible and boundary irreducible 3-manifolds.

\begin{prop}[\textbf{Main result}] \label{main_result}
For every compact irreducible and boundary-irreducible 3-manifold $M$,
$$ c_{GM}(M) = c_{HM}(M),$$
where $ c_{GM}(M)$ (resp. $c_{HM}(M))$ denotes the Gem-Matveev complexity (resp. the modified Heegaard complexity) of $M$.
\end{prop}

The possibility of effectively computing GM-complexity in a direct and algorithmic way from any graph representing $M$ or its associated singular manifold (see Section 5) allows to extend to the boundary case the search for estimations of the Matveev's complexity via $c_{GM}(M)$, as already obtained in the closed case for interesting classes of 3-manifolds:
in particular, in \cite{[C$_5$]}, $GM$-complexity has produced significant improvements in order to estimate Matveev's complexity for two-fold branched coverings of $\mathbb S^3$, for three-fold simple branched coverings of $\mathbb S^3$ and for 3-manifolds obtained by Dehn surgery on framed links in $\mathbb S^3$.

In the last section of the present paper a noteworthy estimation of Matveev's complexity is obtained for all Seifert 3-manifolds with base $\mathbb D^2$ and two exceptional fibers and therefore, as a particular case, of all torus knot complements: see Proposition \ref{complessità grafi parametrici}.

In some particular cases (including for example, three torus knot complements with complexity one and the orientable I-bundle over the Klein bottle, having complexity zero), the obtained estimation turns out to coincide with the exact value of Matveev's complexity (see Corollary \ref{complessità torus knot specifici}).

\bigskip
\bigskip

\section{\hskip -0.7cm . Modified Heegaard complexity}\label{c_HM}

\bigskip

The notion of \textit{modified Heegaard complexity} for compact
orientable 3-manifolds (either with or without boundary) has been
introduced in \cite{[CMV]}, where a comparison with Matveev's
complexity has been discussed. The extension to the non-orientable case has been performed in \cite{[CCM]}. In the present section we recall the main elements of the subject, in the widest possible setting.

Let $\Sigma_g$ be either the closed, connected orientable surface
$T_g$ of genus $g$ (with $g\ge 0$) or the closed, connected
non-orientable surface $U_{2g}$ of genus $2g$ (with $g\ge 1$). So
$\Sigma_g$ is the boundary of a handlebody $\mathbb X_g$ of genus
$g$, $\mathbb X_g$ being the orientable (resp. non orientable)
3-manifold obtained from the $3$-ball $\mathbb D^3$ by adding $g$
orientable 1-handles (resp. $g$ 1-handles, of which at least one
is non-orientable).

A \textit{system of curves} on $\Sigma_g$ is a (possibly empty)
set of simple closed orientation-preserving\footnote{This means
that each curve $\gamma_i$ has an annular regular neighborhood, as
it always happens if $\Sigma_g$ is an orientable surface.} curves
$\mathcal C=\{\gamma_1,\ldots,\gamma_k\}$ on $\Sigma_g$ such that
$\gamma_i \cap \gamma_j = \emptyset$, for $1\le i\neq j\le k$.
Moreover, we denote by
$V(\mathcal C)$ the set of connected
components of the surface obtained by cutting $\Sigma_g$ along the
curves of $\mathcal C$. The system $\mathcal C$ is said to be
\textit{proper} if all elements of $V(\mathcal C)$ have genus
zero, and \textit{reduced} if either $\vert V(\mathcal C)\vert =1$
or no element of $V(\mathcal C)$ has genus 0. Note that a proper reduced system of curves on
$\Sigma_g$ contains exactly $g$ curves.

Let $G(\mathcal C)$ denote the graph which is dual to the one
determined by $\mathcal C$ on $\Sigma_g$; thus, vertices of $G(\mathcal
C)$ correspond to elements of $V(\mathcal C)$ and edges correspond to
curves of $\mathcal C$. Note that loops and multiple edges may arise
in~$G(\mathcal C)$.

A \textit{compression body} $\mathbb Y_g$ of genus $g$ is a 3-manifold with
boundary obtained from $\Sigma_g\times [0,1]$ by attaching a finite
set of 2-handles along a system of  curves (called
\textit{attaching circles}) on $\Sigma_g\times\{0\}$ and filling in
with balls all the spherical boundary components of the resulting
manifold, except $\Sigma_g\times\{1\}$ when $g=0.$
Moreover, $\partial_+ \mathbb Y_g=\Sigma_g\times\{1\}$ is called the
\textit{positive} boundary of $\mathbb Y_g$, while $\partial_- \mathbb Y_g =
\partial \mathbb Y_g-\partial_+ \mathbb Y_g$ is called the \textit{negative}
boundary of $\mathbb Y_g$. A compression body turns out to be a handlebody
if and only if $\partial_- \mathbb Y_g = \emptyset$, i.e., the system of
the attaching circles on $\Sigma_g\times\{0\}$ is proper.
Obviously homeomorphic compression bodies can be obtained via
(infinitely many) non-isotopic systems of attaching circles.

If a system of attaching circles $\mathcal C$ yielding $\mathbb Y_g$
is not reduced, then it contains at least one reduced subsystem of curves  determining the
same compression body $\mathbb Y_g$. Indeed, let  $V^+(\mathcal C)$ be the set of
vertices of $G(\mathcal C)$ corresponding to the components with genus
greater then zero, and $\mathcal A(\mathcal C)$ be the set
consisting of  all the subgraphs $T$ of $G(\mathcal C)$
such that:
\begin{itemize}
\item $T$ contains
all vertices of $G(\mathcal C)$;
\item if $V^+(\mathcal C)=\emptyset,$  $T$ is a (maximal) tree in $G(\mathcal C)$;
\item if $V^+(\mathcal C)\ne \emptyset,$ each connected component of $T$ is a tree
containing exactly one vertex of $V^+(\mathcal C)$.
\end{itemize}
Then, for any  $T \in \mathcal A(\mathcal C)$, the system of
curves obtained by removing from $\mathcal C$ the curves
corresponding to the edges of $T$ is reduced and determines the
same compression body. Note that this operation corresponds to
removing complementary 2- and 3-handles. Moreover, if  $\partial_-
\mathbb Y_g$ is orientable (resp. non-orientable) and has $h$ boundary
components with genus
$^{\partial}g_j$ (resp. $2 \cdot \, ^{\partial}\! g_j$), $1\le j \le h$,
then
$$
\vert E(T)\vert = \vert\mathcal C\vert - g - \max\{0,h-1\} +
\sum_{j=1}^h \ ^{\partial}\! g_j
$$
for each $T\in\mathcal A(\mathcal C)$, where $E(T)$ denotes the edge
set of $T$.

Let $M$ be  a compact connected 3-manifold without spherical boundary components. A
\textit{Heegaard surface} for $M$ is a surface $\Sigma_g$
embedded in $M$ such that $M-\Sigma_g$ consists of two components
$X$  and $Y$ whose closures are (homeomorphic to) a genus $g$
handlebody and  a genus $g$ compression body, respectively.

The triple $(\Sigma_g, X,Y)$ is called a
\textit{(generalized) Heegaard splitting} of genus $g$ of $M$. It is a well
known fact that each compact connected 3-manifold without
spherical boundary components admits a Heegaard splitting.

\begin{rem}
\textup{By Proposition 2.1.5 of \cite{[Ma]}, the complexity of a
manifold is not affected by puncturing it. So, in order to compute
complexity, no loss of generality occurs by assuming that the manifold
has no spherical boundary components.}
\end{rem}

On the other hand, a triple $\mathcal H=(\Sigma_g, \mathcal C',\mathcal C'')$, where  $\mathcal C'$ and
$\mathcal C''$ are two systems of curves on $\Sigma_g$, such that
they intersect transversally and $\mathcal C'$ is proper, uniquely
determines a 3-manifold $M_{\mathcal H}$ corresponding to the (generalized) Heegaard
splitting $(\Sigma_g, X, Y)$, where $X$ and $Y$ are respectively
the handlebody and the compression body whose attaching circles
correspond to the curves in the two systems. Such a triple is called a
\textit{generalized Heegaard diagram}  for $M_{\mathcal H}$.\footnote{In the case of closed 3-manifolds, both systems of curves of a generalized Heegaard diagram $\mathcal H$ are obviously proper; if they are also reduced, $\mathcal H$ is simply a Heegaard diagram in the classical sense (see \cite{[He]}).}

\smallskip

For each generalized Heegaard diagram  $\mathcal H=(\Sigma_g, \mathcal C',\mathcal C'')$, we denote by $\Delta(\mathcal H)$  the graph embedded in $\Sigma_g$ defined by the curves of $\mathcal C'\cup \mathcal C''$, and by $\mathcal R(\mathcal H)$ the set of regions of $\Sigma_g-\Delta(\mathcal H)$. Note that $\Delta(\mathcal H)$ may have connected components which are circles. All vertices not belonging in these components are 4-valent and they are called \textit{singular} vertices.  A diagram $\mathcal H$ is called \textit{reduced Heegaard diagram} if both systems of curves are reduced.
If $\mathcal H$ is non-reduced, any reduction of both its systems of curves yields a reduced Heegaard diagram, which is said to be {\it obtained} from  $\mathcal H$.

The modified complexity of a reduced Heegaard diagram  $\mathcal H'$ is
$$c_{HM}(\mathcal H')= n(\mathcal H')-\max\,\{m(R)\mid R\in\mathcal R(\mathcal H')\},$$ where $n(\mathcal H')$ is the number of singular vertices of $\Delta(\mathcal H')$ and $m(R)$ denotes the number of singular vertices contained in the
region $R$; while the modified complexity of a (non-reduced) generalized Heegaard
diagram $\mathcal H$ is
$$c_{HM}(\mathcal H)=\min\,\{c_{HM}(\mathcal H')\mid \mathcal H' \text{\ reduced \ Heegaard \ diagram \ obtained \ from \ } \mathcal H \}.$$

We define the \textit{modified Heegaard complexity} of a compact
connected 3-manifold $M$ as
$$
c_{HM}(M) = \min\,\{c_{HM}(\mathcal H) \mid \mathcal H  \text{\ generalized \ Heegaard \ diagram \  of \ } M\}.$$
\medskip

The significance of  modified Heegaard complexity consists in its
relation with Matveev's complexity $c(M)$ (see \cite{[CMV]} for the orientable case and \cite{[CCM]} for  the general one):
\begin{prop}
If $M$ is a compact connected 3-manifold, then
$$
c(M) \leqslant c_{HM}(M) .
$$
\end{prop}

\bigskip
\bigskip

\section{\hskip -0.7cm . Basic notions of crystallization theory}
\label{preliminari cGM}

\bigskip

The present section is devoted to briefly review some basic
concepts of the representation theory of PL-manifolds by a particular type of  edge-coloured graphs, called crystallizations.

For general PL-topology and elementary notions about graphs and embeddings, we
refer to \cite{[HW]} and \cite{[Wh]} respectively.

Crystallization theory is a combinatorial tool for representing  compact PL-manifolds, without assumptions about dimension, connectedness, orientability or boundary properties  (see the survey papers \cite{[FGG]}, \cite{[BCG]} and \cite{[BCCGM]}, together with their wide bibliography).
However, since this paper concerns only dimension $3$, we will restrict definitions and results to this dimension, although they mostly hold for the general case; moreover, from now on all manifolds will be assumed to be compact and connected.

Given a pseudocomplex $K$, triangulating a $3$-manifold $M$, a \textit{coloration} on $K$ is a labelling of its vertices by $\Delta_3=\{0,1,2,3\}$, so that:
\begin{itemize}
\item[-] the labelling is injective on each simplex of $K$;
\item[-] each $3$-labelled vertex is internal in $K$.
\end{itemize}

The dual 1-skeleton of $K$ is a (multi)graph $\Gamma=(V(\Gamma),E(\Gamma))$ embedded in
$|K|=M$; we can define a map $\gamma :
E(\Gamma)\to\Delta_3$ in the following way: $\gamma(e)=c\ $ iff
the vertices of the face dual to $e$ are coloured by
$\Delta_3-\{c\}.$
The map $\gamma$ - which is injective on each pair of adjacent edges of
the graph - is called an \textit{edge-coloration} on $\Gamma,$
while the pair $(\Gamma,\gamma)$ is called a \textit{$4$-coloured
graph representing $M$} or simply a \textit{gem} (where ``gem" stands for  \textit{graph encoded
manifold}: see \cite{[Li]}).
In order to avoid long notations, in the following we will often
omit the edge-coloration when it is not necessary, and we will simply
write $\Gamma$ instead of $(\Gamma,\gamma)$.

Obviously, any $3$-manifold $M$ has a gem inducing it: just take the
barycentric subdivision $H^{\prime}$ of any pseudocomplex $H$ triangulating $M$,
label any vertex of $H^{\prime}$ with the dimension of the open simplex containing it
in $H$, and construct the associated $4$-coloured graph as described above.
Conversely, starting from $\Gamma$, we can always
reconstruct $K(\Gamma)=K$ and hence the manifold $M$:
\begin{itemize}
\item[-] {take a $3$-simplex $\sigma(v)$ for each vertex $v\in V(\Gamma)$;}
\item[-] {for each $i \in \Delta_3$ and for each pair $v,w$ of $i$-adjacent vertices of $\Gamma$, identify the faces of $\sigma(v)$ and $\sigma(w)$ opposite to the $i$-coloured vertices, taking care to identify vertices of the same colour.}
\end{itemize}

It is easy to see that, if $M$ is a closed $3$-manifold, any gem of $M$ is a regular graph of
degree $4$, while any gem $\Gamma$ of a $3$-manifold $M$ with non-empty
boundary has a subset of vertices -  called {\it boundary
vertices} - of degree $3$ - which lack in the $3$-coloured edges and correspond to
the tetrahedra of $K(\Gamma)$ having a boundary face.

Given $i,j\in\Delta_3$, $i\neq j$, we denote by $(\Gamma_{i,j},\gamma_{i,j})$
the $2$-coloured graph obtained from $\Gamma$ by deleting all edges which are not $i$- or
$j$-coloured.
The connected components of $\Gamma_{i,j}$ will be called \textit{$\{i,j\}$-residues} of $\Gamma$, and their number will be denoted by $g_{i,j}$.

As a consequence of the above constructions, a bijection is established
between the set of $\{i,j\}$-residues of $\Gamma$
and the set of 1-simplices of $K(\Gamma)$ whose endpoints
are labelled by $\Delta_3-\{i,j\}$.
Moreover, for each $c\in\Delta_3$, the connected components of the
$3$-coloured graph $\Gamma_{\hat c}$ obtained from $\Gamma$ by
deleting all $c$-coloured edges will be called {\it $\hat c$-residues} of $\Gamma$; they are in bijective correspondence with the $c$-labelled vertices of $K(\Gamma)$ and their number will be denoted by $g_{\hat c}.$

Note that, given a $4$-coloured graph $(\Gamma,\gamma)$, $|K(\Gamma)|$ is a $3$-manifold if and only if, for every $c \in \Delta_3$, each connected component $\Xi$ of $\Gamma_{\hat c}$ represents\footnote{By obvious adjustments
of the previous notions, a $2$-dimensional PL-manifold may be encoded by a $3$-coloured graph, which {\it represents} it.} either $\mathbb S^2$ or $\mathbb D^2:$ in the first (resp. latter) case $\Xi$ corresponds to a internal (resp. boundary) vertex of $K(\Gamma).$\footnote{See \cite{[G$_1$]} for definitions and results about $4$-coloured graphs representing polyhedra which fail to be $3$-manifolds. In particular, in Section 5, we will make use also of $4$-coloured graphs representing {\it singular $3$-manifolds}, and hence we will admit a $\hat c$-residue to represent a surface with genus different from zero.}

We will call $(\Gamma,\gamma)$ \textit{contracted} iff  $K(\Gamma)$ has the the minimal number of vertices. This means that:
\begin{itemize}
\item[-]  $\Gamma_{\hat c}$ is connected for each $c\in\Delta_3$, in case $K(\Gamma)$ has either empty or connected boundary;
\item[-] $\Gamma_{\hat 3}$ is connected  and $\Gamma_{\hat c}$ (for each $c \in \{0,1,2\}$) has exactly $h$ connected components,  in case  $K(\Gamma)$ has $h \ge 2$ boundary components.
\end{itemize}

A contracted $4$-coloured graph representing a $3$-manifold $M$ is called a
\textit{crystallization} of $M$.
It is well-known that every $3$-manifold admits a crystallization (see \cite{[FGG]}, together with its references).
Any crystallization (or more generally any gem) $\Gamma$ of $M$ encodes in a combinatorial way the topological properties of $M$. For example, $M$ is orientable iff $\Gamma$ is bipartite.

\bigskip

In order to understand the strong relation existing between gems and (generalized) Heegaard splittings, even in the boundary case, a
particular type of embeddings of a $4$-coloured graph into a surface has to be introduced (see \cite{[G$_3$]} and \cite{[G$_4$]}, where the following definitions and results are given in arbitrary dimension).

Let $(\Gamma,\gamma)$ be a $4$-coloured graph with non-empty boundary (representing a compact 3-manifold); we call {\it extended graph
associated to\/} $(\Gamma,\gamma)$ the $4$-coloured graph
$(\Gamma ^*,\gamma ^*)$ obtained in the following way:
\begin{itemize}
\item[-] {add to $V(\Gamma)$ a set $V^*$ in bijective correspondence with the
set of boundary vertices of $(\Gamma,\gamma)$;}
\item[-] {add to $E(\Gamma)$ the set of all possible $3$-coloured edges
having as endpoints a boundary vertex of $(\Gamma,\gamma)$ and its corresponding vertex in $V^*$.}
\end{itemize}

A {\it regular imbedding\/} of $(\Gamma,\gamma)$ into a surface (with boundary)
$F$ is a cellular imbedding of $(\Gamma ^*,\gamma ^*)$ into $F$, such that:
\begin{itemize}
\item[-]  {the image of a vertex of $\Gamma ^*$ lies in $\partial F$ iff the
vertex belongs to $V^*;$}
\item[-] {the boundary of any region  of the imbedding is either the image
of a cycle of $(\Gamma ^*,\gamma ^*)$ ({\it internal region\/}) or the union of
the image $\alpha$ of a path in $(\Gamma ^*,\gamma ^*)$ and an arc of $\partial
F$, the intersection consisting of the images of two vertices belonging to $V^*$
({\it boundary   region\/});}
\item[-] {there exists a cyclic permutation $\varepsilon =(\varepsilon_0,\varepsilon_1,
\varepsilon_{2},\varepsilon_3=3)$  of $\Delta _3$ such that for each internal region (resp.
boundary region), the edges of its boundary (resp. of $\alpha$) are
alternatively coloured $\varepsilon _i$ and $\varepsilon _{i+1}$, $i\in\mathbb Z_{4}$.}
\end{itemize}

For each cyclic permutation $\varepsilon$ of $\Delta _3$, a regular imbedding of $(\Gamma,\gamma)$ into a surface (with boundary) $F_\varepsilon$ is proved to exists. $F_\varepsilon$  is called the {\it regular surface\/}  associated to $(\Gamma,\gamma)$ and
$\varepsilon;$  it is orientable iff $(\Gamma,\gamma)$ is bipartite.

In the case of $\G$ bipartite  (resp. non-bipartite), the {\it regular genus} $\rho _\varepsilon(\Gamma)$ of $\Gamma$ {\it with respect to $\varepsilon$} is defined as the  classical genus (resp. half the classical genus) of the surface $F_\varepsilon.$

In any dimension, $\rho _\varepsilon(\Gamma)$ may be computed by a suitable combinatorial formula.
In particular, in dimension $3$, we have:
$$ \rho _\varepsilon(\Gamma)=   g_{\varepsilon _0 \varepsilon _{2}}-g_{\hat {\varepsilon_1}}-g_{\hat 3}+1  \ \ \ \ (\text{or, equivalently: } \   \rho _\varepsilon(\Gamma)=   g_{\varepsilon_1 3}-g_{\hat {\varepsilon_0}}-g_{\hat {\varepsilon_{2}}}+1 ) $$

We point out that, for each compact $3$-manifold $M$, the Heegaard genus of $M$ coincides with the minimum value of  $\rho_\varepsilon (\Gamma)$, where $\Gamma$ is a  gem of $M$ and $\varepsilon$ is a cyclic permutation of $\Delta _3$ (see \cite{[G$_2$]} for the closed case, \cite{[CGG]} and \cite{[Cr]}
for the boundary case).\footnote{Actually, $\mathcal G(M)=\min \{\rho_\varepsilon (\Gamma) \ | \ (\Gamma,\gamma) \ \text { gem \ of  } M, \  \varepsilon \ \text { cyclic permutation of } \Delta _n \}$  is an interesting PL-manifold invariant, called {\it regular genus}, which extends to arbitrary dimension the classical notion of  Heegaard genus of a $3$-manifold (see \cite{[G$_2$]} for details).
By means of regular genus, important classification results have been obtained within crystallization theory: see, for example, \cite{[C$_1$]} and \cite{[CM]} for dimension four, \cite{[CG]},   \cite{[C$_2$]}  and \cite{[C$_3$]} for dimension five.}

\bigskip
\bigskip

\section{GM-complexity for compact 3-manifolds}\label{GM compact}

Let $M$ be a compact 3-manifold and let $\G$ be a gem representing $M$.
Since the theory concerning closed 3-manifolds has already been developed (see \cite{[C$_4$]}, \cite {[C$_5$]} and \cite{[CC$_1$]}),  in the following we will suppose  $\partial M\neq\emptyset$; let $h \ge 1$ the number of connected components of $\partial M$.

If $\alpha\in\Delta_2$ is an arbitrarily fixed colour, we set $\Delta_3=\{\alpha,\beta,\beta^\prime,3\}$ and consider the cyclic permutation $\varepsilon= (\beta,\alpha,\beta^\prime,3)$ of $\Delta_3$.

We denote by $K_{\alpha 3}$  (resp. $K_{\beta\beta^\prime}$) the 1-dimensional subcomplex of $K=K(\Gamma)$ generated by the
$\alpha-$ and $3-$ (resp. $\beta-$ and $\beta^\prime-$) labelled vertices. Let $H_{\alpha}$ be the largest 2-dimensional subcomplex of the first barycentric subdivision of $K$ not intersecting the subdivisions of $K_{\alpha 3}$ and $K_{\beta\beta^\prime}$, and set $F^{(\alpha)}=|H_{\alpha}|$.

The surface (with boundary) $F^{(\alpha)}$ splits $K(\Gamma)$ into two polyhedra $\mathcal A_{\alpha 3}$  and  $\mathcal A_{\beta\beta^{\prime}}$
(see Figure 1). It is not difficult to check - by constructions described in  \cite{[G$_3$]} and \cite{[G$_4$]} - that $F^{(\alpha)}$ is the surface into which $\G$ regularly embeds with respect to the chosen permutation $\varepsilon$.

\medskip
\centerline{\scalebox{1.2}{\includegraphics{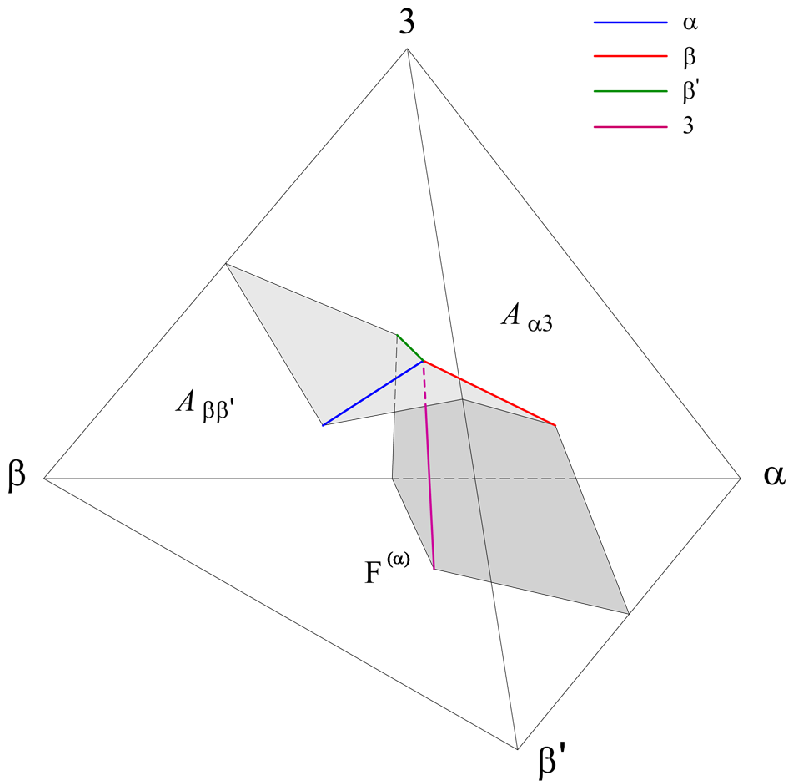}}}
\smallskip
\centerline{\footnotesize{Figure 1: the local splitting of $K(\Gamma)$ into $\mathcal A_{\alpha 3}$  and $\mathcal A_{\beta\beta^{\prime}}$}}

\bigskip

$\mathcal A_{\alpha 3}$ (resp. $\mathcal A_{\beta\beta^\prime}$) is a handlebody, since it collapses to the graph $K_{\alpha 3}$ (resp. $K_{\beta\beta^\prime}$).

Furthermore, note that $\mathcal A_{\alpha 3}\cap\mathcal A_{\beta\beta^\prime}=\partial\mathcal A_{\alpha 3}\cap\partial\mathcal A_{\beta\beta^\prime}=F^{(\alpha)}$.

Since $\partial M\neq\emptyset$, the two handlebodies do not intersect on their whole boundaries; in particular $\partial\mathcal A_{\alpha 3}\cap\partial M=\bigcup_{i=1}^{t}\mathbb D_i$, where the $\mathbb D_i$'s are 2-disks, such that $\partial\mathbb D_i = lk(P_i,(\partial K)^\prime)$ for each $i=1, \dots, t$, where $P_1, P_2, \dots, P_t$ ($t \ge h \ge 1$) are the \ $\alpha-$labelled vertices of $\partial K$ and  $(\partial K)^\prime$ denotes the first barycentric subdivision of $\partial K$.
Moreover, $\bigcup_{i=1}^{t}\partial\mathbb D_i=\partial F^{(\alpha)}$.

Let $\Sigma^{(\alpha)}$ be the closed surface obtained by adding the disks $\mathbb D_1,\ldots,\mathbb D_{t}$
to $F^{(\alpha)}$ and let $C = \Sigma^{(\alpha)}\times [-1,1]$ be a collar of $\Sigma^{(\alpha)}$ in $\mathcal A_{\alpha 3}$.

We set $C^- = \Sigma^{(\alpha)} \times [-1,0],$ \ $C^+ = \Sigma^{(\alpha)} \times [0,1]$ and $S_j = \Sigma^{(\alpha)} \times\{j\}$ for each $j \in \{-1, 0, 1\}$ (see Figure 2).

Let us define the pseudocomplexes:  $X_\alpha=\overline{\mathcal A_{\alpha 3}\setminus C^-}$ and $Y_\alpha=\mathcal A_{\beta\beta^\prime}\cup C^-$.

$X_\alpha$  is obtained from $C^+$ by attaching 2-handles to $S_1$  along the  $\{\beta,\beta^\prime\}$-coloured cycles of  $\G$ (dual to the 1-simplices of $K_{\alpha 3}$).
Hence, $X_\alpha$ is a compression body, whose system of attaching circles is the set $\mathbf {v_\alpha}$ of $\{\beta,\beta^\prime\}$-coloured cycles of $\G$ (dual to edges of $K$ with $\alpha$- and $3$-coloured endpoints).  $\mathbf {v_\alpha}$ is a proper, not reduced system of curves on the surface $S_0\cong \Sigma^{(\alpha)}$: in fact, all connected components of $\Sigma^{(\alpha)}\setminus \mathbf {v_\alpha}$ have genus zero, since they are the boundaries of suitable neighborhoods either of the $\alpha$-labelled vertices of $\partial K$ (hence 2-disks)  or  of the internal $\alpha$- and $3$-labelled
vertices (hence 2-spheres with holes).
As a consequence, $X_\alpha$  turns out to be a  handlebody of genus $g(\Sigma^{(\alpha)})=g(F^{(\alpha)})=g_{\beta\beta^\prime}-g_{\hat\alpha}-g_{\hat 3}+1$.
The positive boundary of $X_\alpha$ is $S_0$ (while its negative boundary is obviously empty).

\medskip
\centerline{\scalebox{0.5}{\includegraphics{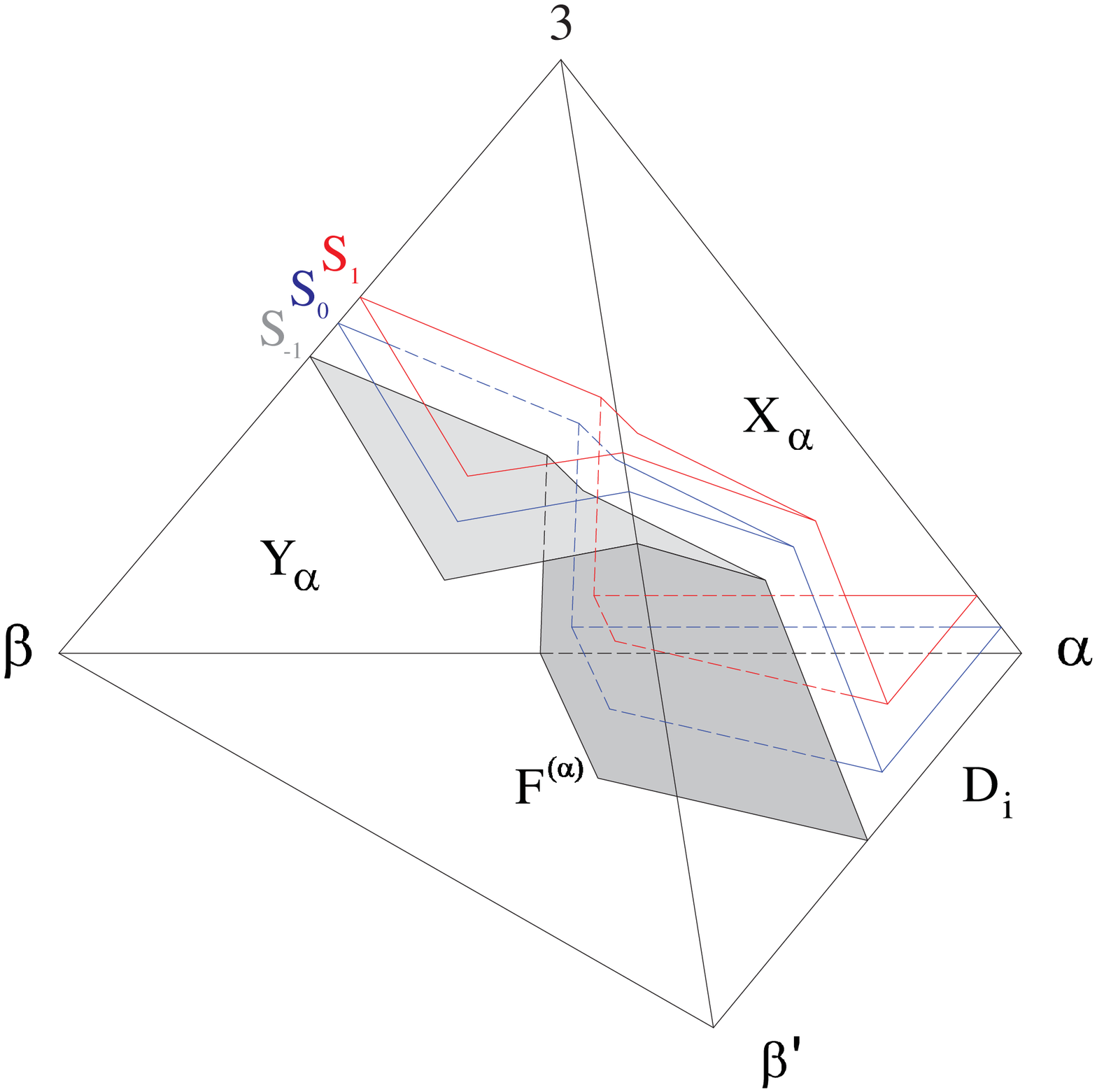}}}
\smallskip
\centerline{\footnotesize{Figure 2: the ``double" collar of $\Sigma^{(\alpha)}$ in $\mathcal A_{\alpha 3}$}}

\bigskip

On the other hand, $Y_\alpha$ is a compression body, obtained from $C^-$ by attaching 2-handles to $\Sigma^{(\alpha)} = S_{-1}$ along  the  $\{\alpha,3\}$-coloured cycles of $\G$ (dual to 1-simplices of  $K_{\beta\beta^\prime}$ not belonging to $\partial K$).
The set $\mathbf {w_\alpha}$  of $\{\alpha,3\}$-coloured cycles of $\G$ is a not proper (if $\partial M\neq\emptyset$ with at least one non-spherical component) and generally not reduced (unless $\G$ is a crystallization) system of curves on $S_0\cong \Sigma^{(\alpha)}$.

Therefore, the triple $(\Sigma^{(\alpha)},X_\alpha,Y_\alpha)$  (resp. $\mathcal H_\alpha=(\Sigma^{(\alpha)},\mathbf {v_\alpha},\mathbf {w_\alpha})$) \ is a (generalized) Heegaard splitting (resp. generalized Heegaard diagram) of $M$.

\smallskip

In order to get reduced systems of curves, i.e. to obtain a reduced Heegaard diagram for $M$ starting from $ \mathcal H_\alpha$, let us first denote by $G_{\beta\beta^\prime}$ the graph obtained from $K_{\beta\beta^\prime}$ by contracting to one point $Q_i$ (for each $i=1,\ldots,h$) the vertices of $K_{\beta\beta^\prime}$ belonging to the  $i$-th component of $\partial K$.

We consider the following sets of curves on $\Sigma^{(\alpha)}$:
\begin{itemize}
\item [(a)] $\mathcal D$, which is a set of  $\{\beta, \beta^\prime\}$- coloured cycles dual to a maximal tree of  $K_{\alpha 3}.$

\item [(b)]  $\mathcal D^\prime $, which is a set of $\{\alpha, 3\}$-coloured cycles dual to the edges of a subgraph $\bar G$ of  $G_{\beta\beta^\prime}$ such that $\bar G$ is a union of trees containing all vertices of  $G_{\beta\beta^\prime}$ and, for each $i, j,$ \ $i \ne j$, the vertices $Q_i$ and  $Q_j$ belong to different connected components of $\bar G.$
\end{itemize}

\begin{prop}\label{reduced Heegaard} $\mathcal H_\alpha(\mathcal D,\mathcal D^\prime)=(\Sigma^{(\alpha)},\mathbf {v_\alpha}\setminus\mathcal D,\mathbf {w_\alpha}\setminus\mathcal D^\prime)$ \ \ is a reduced Heegaard diagram of $M$.\end{prop}

\dimo The graph $K_{\alpha 3}$ is isomorphic to the graph $G(\mathbf {v_\alpha})$, which is the dual 1-skeleton of the cellular decomposition induced by the curves of $\mathbf {v_\alpha}$ on $S_0\cong \Sigma^{(\alpha)}$.
Therefore, by removing from $\mathbf {v_\alpha}$ the curves of $\mathcal D$, we get a reduced system of attaching curves for the 2-handles on $S_1$, yielding the handlebody $X_\alpha$.

On the other hand, $G_{\beta\beta^\prime}$ is isomorphic to the graph $G(\mathbf {w_\alpha})$, which is the dual 1-skeleton of the cellular decomposition induced by the curves of $\mathbf {w_\alpha}$ on $S_0\cong \Sigma^{(\alpha)}$.
Hence, removing the curves of $\mathcal D^\prime$ from $\mathbf {w_\alpha}$ does not affect $Y_\alpha$, since it corresponds to the cancellation of complementary 2- and 3-handles.
\qed

\begin{defn} \label{def c_GM} \emph{Given a gem $\G,$ with the above notations, we define the \textit{Gem-Matveev complexity} ({\it GM-complexity} for short) of $\G$ as
\begin{eqnarray*}
c_{GM}(\G)= \min\{c_{HM}(\mathcal H_    \alpha(\mathcal D,\mathcal D^\prime))\ |\ \alpha\in\Delta_2, \mathcal D\subset\G_{\alpha 3},\mathcal D^\prime\subset\G_{\beta \beta^\prime} \ \text{\ satisfying (a) e (b)}\}
\end{eqnarray*}} \end{defn}

\bigskip
We point out that, although the above definition of GM-complexity is given through the associated generalized Heegaard diagrams, nevertheless it can be computed algorithmically and directly from the graph, as it happens in the closed case (\cite{[C$_4$]}).

\begin{rem}\label{cryst} \emph{If $\G$ is crystallization, then a reduced system of curves for the handlebody  $X_\alpha$ is obtained by considering the connected components $\Xi_i$ ($i=1,\ldots, h$)\footnote{Recall that, if $\G$ is a crystallization of $M$ and $\alpha \in \Delta_2,$  $\G_{\hat \alpha}$ consists of exactly $h$ connected 3-coloured graphs with non-empty boundary, where $h$ is the number of components of $\partial M$.} of the subgraph $\G_{\hat \alpha}$ containing boundary vertices and by removing from $\mathbf {v_\alpha}$ a $\{\beta,\beta^\prime\}$-coloured cycle belonging to $\Xi_i$, for each $i=1,\ldots,h$.
The new system has exactly $g_{\beta\beta^\prime}- h =g(\Sigma^{(\alpha)})$ curves.
On the other hand, in this case,
 $\mathbf {w_\alpha}$ is an already reduced system.}\end{rem}

\begin{defn}\emph{Given a compact 3-manifold $M$, the \textit{Gem-Matveev complexity} (GM-complexity for short) of $M$ is defined as
$$c_{GM}(M)= \min \{ c_{GM}(\G)\ |\ \G\text{\ \ gem of  } M\}.$$}\end{defn}

\medskip

\begin{rem}\label{spherical}\emph{Let $\bar M$ be a compact 3-manifold whose boundary has only spherical components and let  $M$ be the closed manifold obtained by capping off with a 3-ball each boundary component of $\bar M$. Given a gem $\bar\G$ of $\bar M$, let us suppose $c_{GM}(\bar\G)=c_{HM}(\mathcal H_\alpha)$, where $\mathcal H_\alpha$ is the generalized Heegaard diagram associated to $\bar\G$ and $\alpha\in \Delta_2$.
Let us fix a colour $i\in\Delta_2\setminus\{\alpha\}$ and consider the 4-coloured graph without boundary $\G$ obtained by joining two boundary vertices $u$ and $u^\prime$ of $\bar\G$ by a $3$-coloured edge iff there exists a $\{i,3\}$-coloured path in $\bar\G$ having $u$ and $u^\prime$ as endpoints. This is a general construction introduced in \cite{[G$_2$]}: starting from a gem of a compact 3-manifold, it produces a gem of the (possibly singular) manifold obtained by capping off each boundary component by a cone.
A direct computation of the Euler characteristic shows that the choice of $i\neq\alpha$ guarantees that  $\G$ regularly embeds into the closed surface $\Sigma^{(\alpha)}$, too.
Moreover, it is easy to see that $\G$ admits a reduced Heegaard diagram, which can be obtained from $\mathcal H_\alpha$ and whose HM-complexity is $c_{HM}(\mathcal H_\alpha)$. Hence $c_{GM}(\G)=c_{GM}(\bar G)$.
Therefore, our definition of GM-complexity for 3-manifolds with non-empty boundary is consistent with the analogous definition for the closed case.
Conversely, it is easy to see that, by puncturing a closed manifold, we do not change its GM-complexity, as it happens with Matveev's complexity (see Remark 1).}\end{rem}

\bigskip

Since the GM-complexity of a gem is, by definition, the HM-complexity of one of its associated generalized Heegaard diagrams, directly from Proposition \ref{reduced Heegaard} we have:

\medskip

\begin{lemma}\label{lemma1} For each compact 3-manifold $M$,
$$c_{HM}(M) \, \le \, c_{GM}(M)$$
\end{lemma}

\vskip-0.8truecm

\ \qed

\medskip

\noindent {\bf Example 1:} Let $M= T_1 \times I$ be the
compact orientable 3-manifold obtained by making the product
between the closed orientable genus one surface $T_1$ (i.e., the
torus) and  the unit interval $I$. A crystallization $\Gamma$ of
$M$ is given in \cite[Fig. 2]{[G$_0$]}. Figure 3(a) shows the regular
embedding of $\Gamma$ into the orientable surface
$F_\varepsilon,$ where $\varepsilon =(0,1,2,3);$  note that
$F_\varepsilon$ has  genus two  and two boundary components, corresponding to the
two connected components of $\partial M.$

\medskip
\centerline{\scalebox{0.9}{\includegraphics{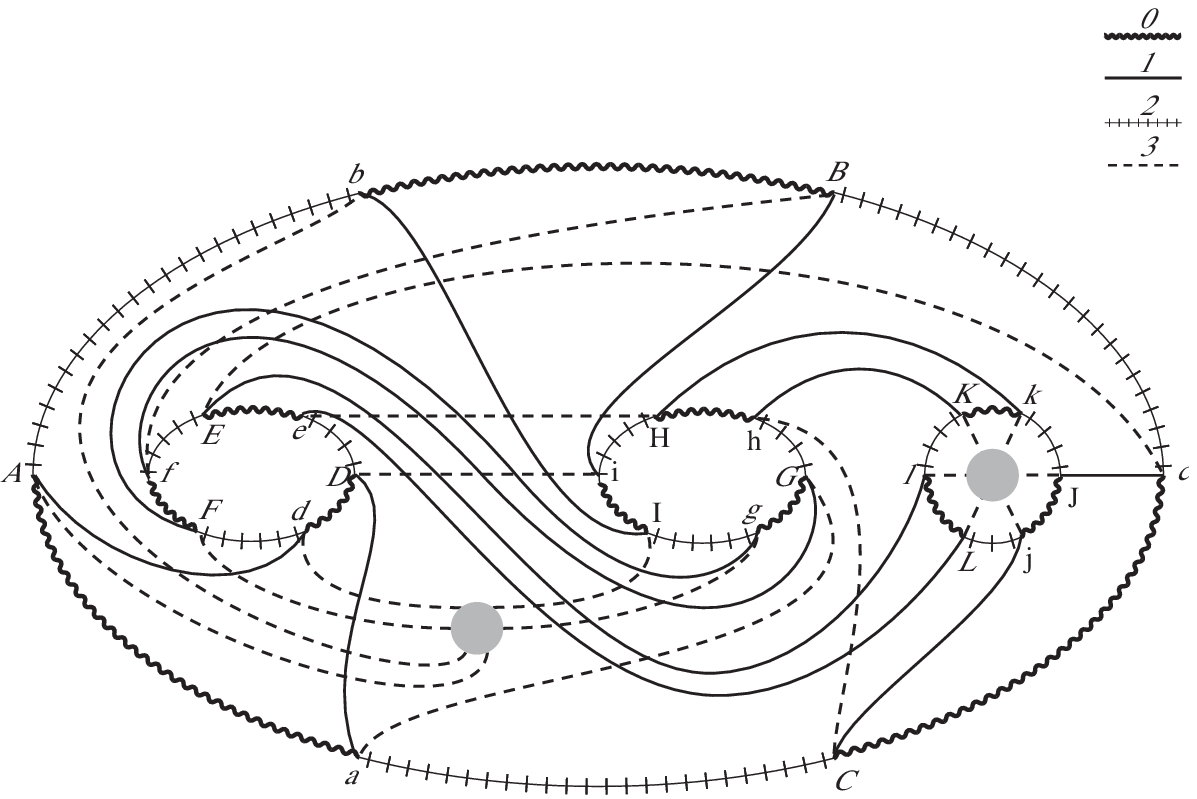}}}
\smallskip
\centerline{\footnotesize{Figure 3(a): a crystallization of $T_1 \times I$}}

\bigskip

By Definition \ref{def c_GM} and Remark \ref{cryst}, $c_{GM}(\Gamma)$
may be computed by considering, on the splitting surface
$F_\varepsilon,$ the set $\bf v$ of all $\{0,2\}$-coloured cycles
of $\Gamma,$ but one for each connected component of $\G_{\hat
3}$, and the set $\bf w$ of all $\{1,3\}$-coloured cycles of
$\Gamma$. This is equivalent to considering the reduced Heegaard
diagram $\mathcal H_1(\mathcal D,\mathcal D^\prime)$ of
Proposition \ref{reduced Heegaard}, where  $\mathcal D$ consists
of the $\{0,2\}$-coloured cycles having as vertex sets
$\{a,A,b,B,c,C\}$ and  $\{j,J,k,K,l,L\}$ (which are dual to a
maximal tree of  $K_{13}$) and $\mathcal D^\prime = \emptyset.$
Figure 3(b) shows the intersection points of the curves of
this reduced diagram and  indicates a region whose boundary contains all the intersection points; hence, $c_{GM}(\Gamma) =0$ follows
(thus proving the sharpness of the upper bound $c_{GM}(T_1 \times I)=c_{HM}(T_1 \times I)$ with respect to
$c(T_1 \times I)$).

\smallskip
\centerline{\scalebox{1}{\includegraphics{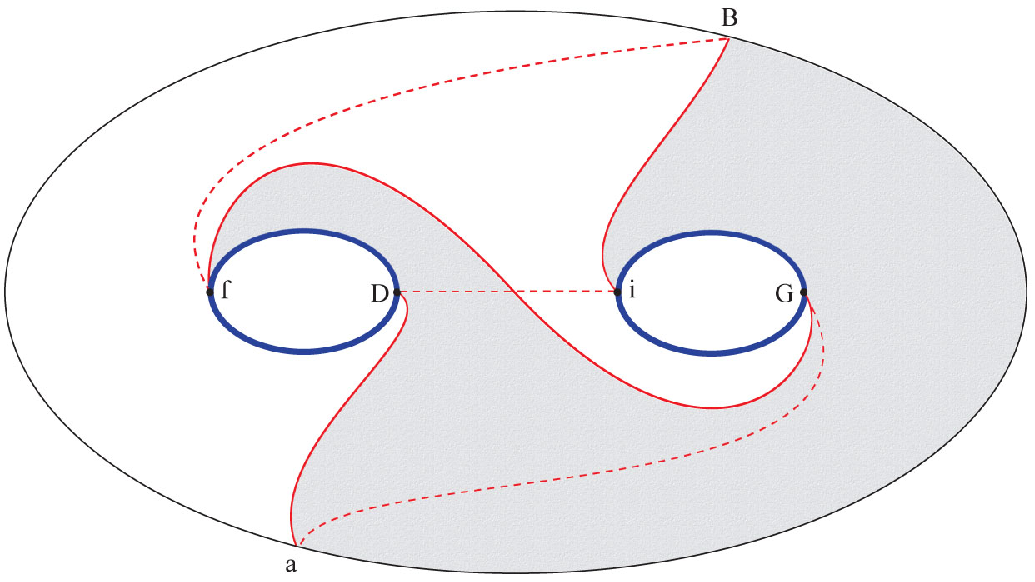}}}
\smallskip
\centerline{\footnotesize{Figure 3(b): a reduced Heegaard diagram of $T_1 \times I$}}

\bigskip
\bigskip

\section{\hskip -0.7cm . GM-complexity via singular 3-manifolds}

\bigskip

In this section, we will extend to singular 3-manifolds some notions of sections 3 and 4; they will be a useful tool in order to prove our main result (Section 6) and to obtain an estimation of Matveev's complexity for all torus knot complements (Section 7).

\bigskip

Let $\G$ be a {\it regular} 4-coloured graph.
The construction of the pseudocomplex $K=K(\G)$ described in section \ref{preliminari cGM}
- which can be obviously performed even if $\G$ does not encode a (closed) 3-manifold -
establishes a bijective correspondence between the $\hat c$-residues
($c\in\Delta_3$) of $\G$ and suitable neighborhoods of the vertices of $K$:
in fact, any  $\hat c$-residue, dual to a vertex $v_c\in V(K)$, represents the surface $|lk(v_c, K^{\prime})|$
(see  \cite{[FGG]}).
Moreover, the cone $v_c*lk(v_c, K^{\prime})$ is a regular neighborhood of $v_c$ in $K$, called the \textit{disjoint star} $std(v_c,K)$ of $v_c$ in $K$.

Hence, for any (regular) $\G$, $|K|$ turns out to be a singular 3-manifold, whose singular points are precisely the vertices of $K$ which correspond to 3-residues of $\G$ not representing 2-spheres (see \cite{[G$_1$]}).

In the following we will always suppose that all the singular vertices of $K$  are labelled by the same colour (which is called the \textit{singular colour} of both $\G$ and $K$).
For sake of simplicity,  we will restrict ourselves to the class $G^{(0)}$ of (regular) 4-coloured graphs $\G$ with singular colour $0$, i.e. such that all $\hat c$-residues of $\G$ represent 2-spheres except when $c=0$.

\medskip

For each $\G\in G^{(0)}$, let us denote by $N= |K(\G)|$ the associated singular 3-manifold  and by $\mathcal S(K)$ the (non-empty) set of singular vertices of $K=K(\G)$.

Note that a compact 3-manifold $M= \check N$ (having non-empty boundary and no spherical boundary component) is obtained  by removing from $N$ small open neighborhoods of its singular points; hence, we will say - by little abuse - that $\G\in G^{(0)}$ {\it represents} $M$, too.  Conversely, the singular manifold  $N=\widehat{M}$ is (uniquely) obtained from $M$ by capping off each component of $\partial M$ by a cone.

\medskip

Let now $\Delta_3=\{0,\alpha,\alpha^\prime,\beta\}$,  let $K_{\alpha\alpha^\prime},\ K_{0\beta}$ and $H_{\alpha}$ be the analogous of the subcomplexes defined in section \ref{GM compact} (in order to recover those definitions, we have simply to set $\alpha^\prime=3$ and $0=\beta^\prime$) and let  $F^{(\alpha)}$ be the (closed) surface triangulated by the 2-complex $H_{\alpha}$.

Furthermore, we denote by $\mathcal A_{\alpha\alpha^\prime}$  and $\mathcal A_{0\beta}$ the two subpolyhedra of $K$ such that  $\mathcal A_{\alpha\alpha^\prime}\cap\mathcal A_{0\beta}=F^{(\alpha)}$ and, obviously, $K_{\alpha\alpha^\prime}\subset\mathcal A_{\alpha\alpha^\prime},\ K_{0\beta}\subset\mathcal A_{0\beta}$.
It is easy to see that  $\mathcal A_{\alpha\alpha^\prime}\cup\mathcal
A_{0\beta}=K$ and  $\mathcal A_{\alpha\alpha^\prime}\cap\mathcal
A_{0\beta}=\partial\mathcal A_{\alpha\alpha^\prime}\cap\partial\mathcal
A_{0\beta}$.
As in section \ref{GM compact}, we point out that $F^{(\alpha)}$ is the surface into which $\G$ regularly embeds with respect to the cyclic permutation $\varepsilon=(\beta,\alpha,0,\alpha^\prime)$.

\bigskip

Note that $X_\alpha=\mathcal A_{\alpha\alpha^\prime}$ is a handlebody, since both the $\alpha$-coloured and the $\alpha^\prime$-coloured vertices of $K$ are not singular.
The set $\mathbf {v_{\alpha}}$ of $\{\beta,0\}$-coloured cycles of $\G$ is a proper, not reduced system of curves on $F^{(\alpha)}$ defining $X_\alpha$ as a compression body.
In order to get a reduced system from $\mathbf {v_{\alpha}}$, we have, exactly as in section \ref{GM compact}, to remove a set $\mathcal D$ of cycles which are dual to the edges of a maximal tree of $K_{\alpha\alpha^\prime}$.

Now, for each $v_0\in\mathcal S(K)$, let us consider a collar $C_{v_0}$ of $lk(v_0,K')$ in the disjoint star of $v_0$ in $K$; let us define the following subset of $\mathcal A_{0\beta}$:
$$\mathcal U = \bigcup_{v_0\in\mathcal S(K)}(std(v_0,K)\setminus C_{v_0})$$

\noindent and  $Y_{\alpha}=\overline{\mathcal A_{0\beta}\setminus\mathcal U}$.
Note that $|X_\alpha\cup Y_\alpha|\cong M$.

$Y_{\alpha}$ is a compression body, since it is obtained from a suitable collar of $F^{(\alpha)}$ in  $\mathcal A_{0\beta}$ by adding 2-handles which have the $\{\alpha,\alpha^\prime\}$-coloured cycles of $\G$ as cores.
Hence, the set  $\mathbf {w_\alpha}$  of such cycles is a non-proper (if $K$ has at least one singular vertex) and (generally) non-reduced system of curves on $F^{(\alpha)}$ representing $Y_{\alpha}$.
In order to reduce $\mathbf {w_\alpha}$, we have to proceed as described in section \ref{c_HM}, by removing a subset $\mathcal D^\prime $ whose curves are dual to the edges of a subgraph $T$ of $G(\mathbf {w_\alpha})$ satisfying the conditions listed in the cited section.

As a consequence, the triple $\mathcal H_{\alpha}(\mathcal D,\mathcal D^\prime)=(F^{(\alpha)},\mathbf {v_\alpha}\setminus\mathcal D,\mathbf {w_\alpha}\setminus\mathcal D^\prime)$ \ is a reduced Heegaard diagram of $M$.

\smallskip

It is thus possible to define an upper bound for the complexity of $M$ through the graphs representing $\widehat{M}$ and having only one singular colour.
More precisely, for each $\G\in G^{(0)}$, we define the {\it GM-complexity of $\G$} (denoted as usual by $c_{GM}(\G)$) as the minimum complexity of the diagrams $\mathcal H_\alpha(\mathcal D,\mathcal D^\prime)$, where $\alpha\in\{1,2,3\}$ and $\mathcal D,$ $\mathcal D^\prime$ are as described above.

Then, we set
$$\hat c_{GM}(M)=\min\{c_{GM}(\G)\ |\ \G\in G^{(0)},\
|K(\G)|\cong\widehat{M}\}.$$

\medskip

By definition itself, we have:
\begin{prop} \label{complessità_singolare} For each compact 3-manifold $M$, \ $c(M) \le \hat c_{GM}(M)$.
\ \qed  \end{prop}

\smallskip

At the end of the next section, as a consequence of the proof of the main result, we will establish the coincidence between GM-complexity computed on edge-coloured graphs representing 3-manifolds and GM-complexity computed on edge-coloured graphs representing the associated singular manifold
(Proposition \ref{uguaglianza complessità varieta' singolari}).

\bigskip

The following Lemma shows how to obtain a gem of the compact 3-manifold $M$ directly from a graph of $G^{(0)}$ representing the singular 3-manifold $\widehat M$; it will be particularly useful in the next section.

\begin{lemma}\label{bis-tris}  \emph{\cite[Lemma 3]{[CGG]}}
If $\G\in G^{(0)}$ represents $\widehat M,$ then a gem $\bar\G$ of  $M$ exists, such that $\rho_\varepsilon(\bar\G)=\rho_\varepsilon(\G)$, where $\varepsilon=(0,1,2,3)$.\end{lemma}
\vskip-0.8truecm
\ \qed

In order to give a precise description of the above graph $\bar\G,$  we have to introduce two particular types of subdivision of a coloured complex: {\it bisection} and {\it trisection} (\cite{[G$_1$]}).
The local effect of  bisection (resp. trisection) of type $(i,j)$ ($i,j\in\Delta_3$) around an $i$-labelled vertex (resp. on a 1-simplex with $i-$ and $j$-labelled endpoints) of a coloured complex is shown in Figure 4, arrow $b(i,j)$ (resp. Figure 6, arrow $t(i,j)$).

\medskip

\centerline{\scalebox{0.75}{\includegraphics{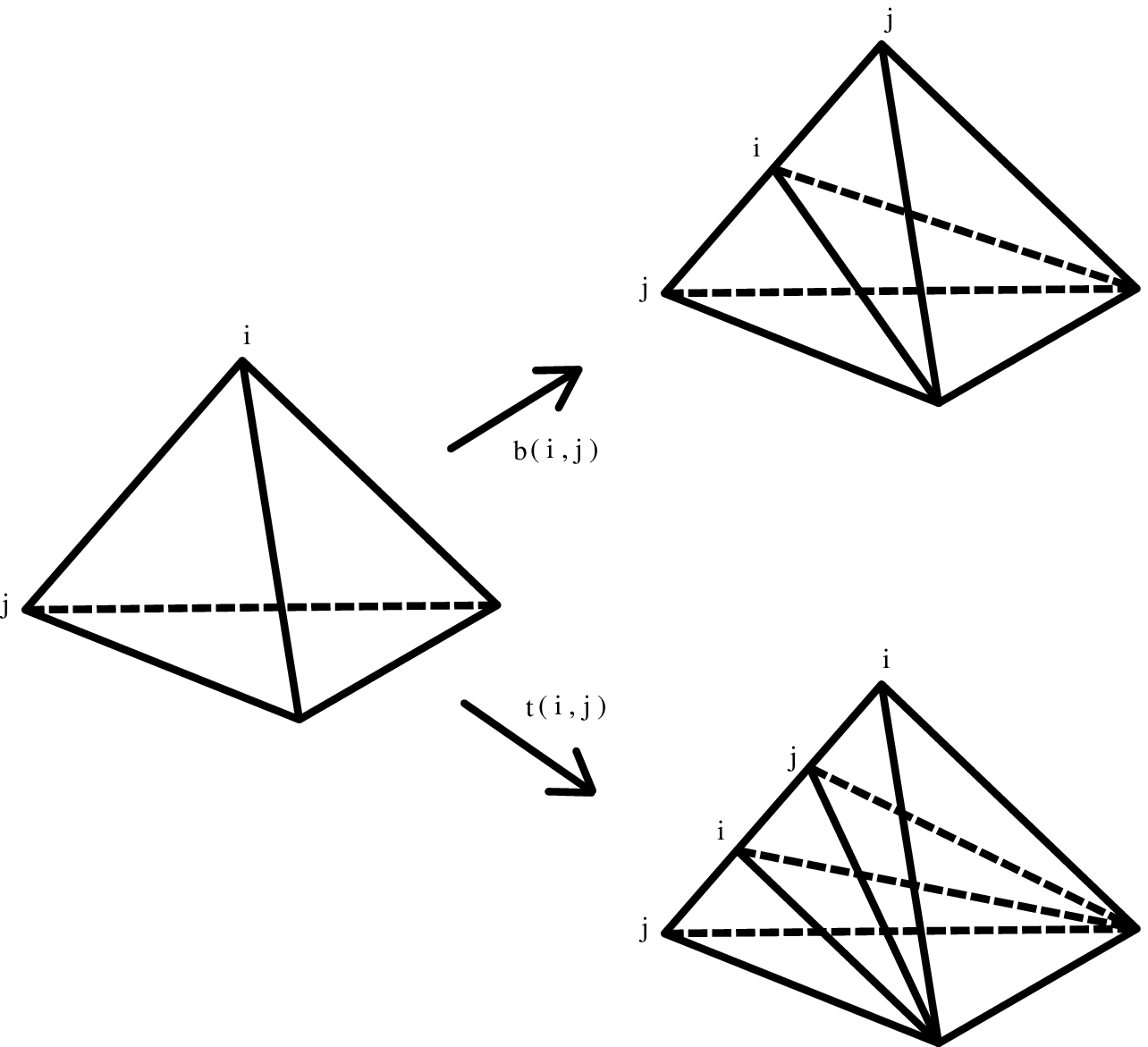}}}
\smallskip

\centerline{\footnotesize{Figure 4: local effect of a bisection (resp. trisection) of type $(i,j)$}}
\bigskip

The graph $\bar\G$ is obtained in the following way:
\begin{itemize}
\item [-] we perform subsequently a trisection of type  $(0,2)$, a bisection of type $(0,3)$ and a trisection of type $(3,1)$ on each edge of $K(\G)$ having a   (0-labelled) singular vertex as endpoint.
    The new coloured complex $\tilde K$ is a coloured subdivision  of $K(\G)$, such that for each vertex $u$ of $\tilde K,$ we have that $|st(u,\tilde K)|$  is the cone over $|lk(u, \tilde K)|$.
  Moreover, the singular vertices of $\tilde K,$  which are exactly those of $K(\G)$, are 3-labelled.
\item [-] we consider the coloured complex $\bar K$ obtained from  $\tilde K$  by deleting  $st(u, \tilde K)\setminus lk(u, \tilde K)$  for each singular  vertex $u$ of
    $\tilde K;$  then,  $\bar\G = \G(\bar K)$.
\end{itemize}

$\bar\G$ turns out to be a gem of $M$, since the removal of open neighborhoods of the singular points of $\widehat M$ yields again $M$.

\bigskip
\bigskip
\section{Proof of the main result}

In this section we will prove that $GM$-complexity and $HM$-complexity coincide for each compact
irreducible and boundary-irreducible 3-manifold, thus extending the analogous result stated for the closed case in \cite{[CCM]}.
Hence, $GM$-complexity and $HM$-complexity  turn out to be useful different tools to compute the same upper bound for Matveev's complexity.

\bigskip

In order to prove our main result, by Lemma \ref{lemma1}, it is sufficient to prove that $c_{GM}(M)\le c_{HM}(M)$, i.e. that each compact irreducible and boundary-irreducible 3-manifold $M$ admits a gem whose GM-complexity is equal to $c_{HM}(M)$.

Let us first observe that we can take into consideration only compact manifolds without spherical boundary components, since puncturing does affect neither GM- nor HM-complexity and the case of manifolds whose boundary has only spherical components is essentially that of closed manifolds (see Remark \ref{spherical}), for which the equality of GM- and HM-complexity has already been established.

Therefore, from now on, we will consider only compact $3$-manifolds, having non-empty boundary and no spherical boundary component.

\bigskip

\centerline{\scalebox{1}{\includegraphics{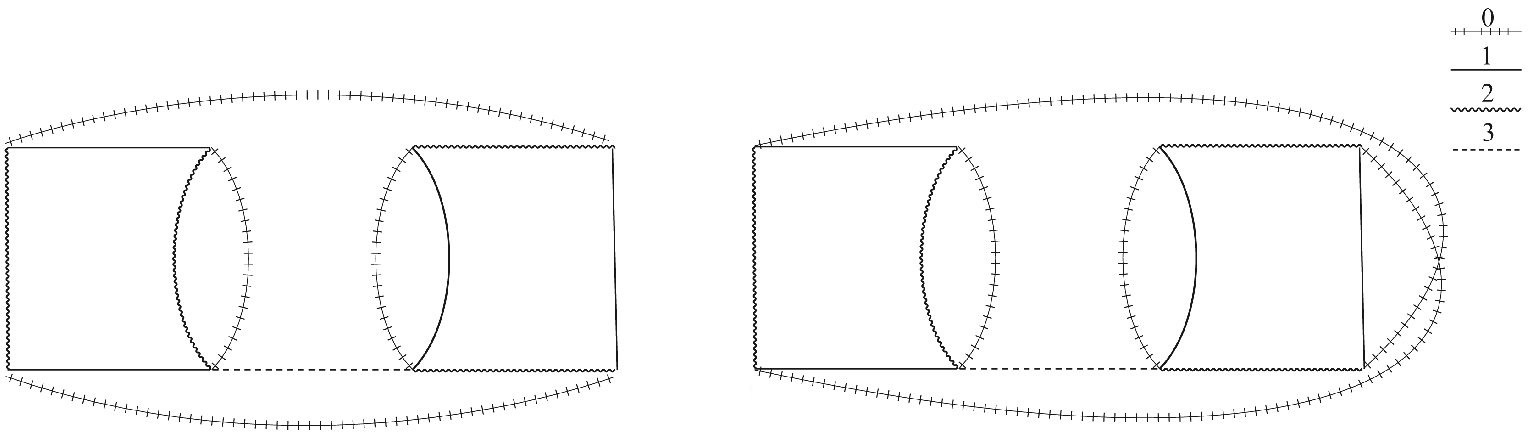}}}
\smallskip
\centerline{\footnotesize{Figure 5: a gem of the orientable (resp. non-orientable) genus $1$ handlebody}}

\bigskip

Let $\mathcal H=(\Sigma_g,\mathbf {v},\mathbf {w})$  be a reduced Heegaard diagram of the compact $3$-manifold $M$, and suppose that $c_{HM}(\mathcal H)=c_{HM}(M)$.

We point out that, if $\mathbf {w}=\emptyset,$ then $M$ is a handlebody and a crystallization $\G$ of $M$ exists such that  $c_{GM}(\G)=c_{HM}(M)=0$
(see Figure 5 for the cases of genus $1$ handlebodies\footnote{Suitable crystallizations of handlebodies with higher genus may easily be obtained by graph connected sum on boundary vertices: see \cite{[FGG]} for details.}); so, from now on, we will suppose $\mathbf {w}\neq\emptyset$.

\bigskip

Furthermore, an Heegaard diagram $\mathcal H=(\Sigma_g,\mathbf v,\mathbf w)$  is called \textit{connected} iff the graph defined on $\Sigma_g$ by the curves of $\mathbf v$ and $\mathbf w$ is connected, i.e. the connected components of $\Sigma_g\setminus (\mathbf v\cup \mathbf w)$ are 2-disks.

Note that if $\mathcal H$ is connected, then it satisfies the following condition:

\medskip
(*)\ \ each curve of $\mathbf v$ intersects at least one curve of $\mathbf w$ and, conversely, each curve of $\mathbf w$ intersects at least one curve of $\mathbf v$.

\bigskip

Obviously, the converse is not always true. In the following we will always refer to the above condition as ``condition (*)".

\bigskip

\begin{lemma}\label{lemma2} Let $\mathcal H$ be a not connected (reduced) Heegaard diagram of a compact 3-manifold $M$, such that $c_{HM}(\mathcal H)=
c_{HM}(M)$. If  $\mathcal H_1,$  $\mathcal H_2,$ $\dots,$ $\mathcal H_k$ ($k \ge 2$) are the connected diagrams into which $\mathcal H$ splits,   then $\sum_{i=1}^k c_{HM}(\mathcal H_i) \le c_{HM}(\mathcal H).$
\end{lemma}

\dimo
Obviously, if $\mathcal H= (\Sigma_g, \mathbf v,\mathbf w),$ then $\mathcal H_i= (\Sigma_{g_i}, \mathbf{v_i},\mathbf{ w_i})$ is such that $\sum_{i=1}^k g_i  =g,$ $\Sigma_{g}= \#_{i=1}^k  \Sigma_{g_i}$ and $\cup_{i=1}^k  \mathbf{ v_i} = \mathbf v$ (resp. $\cup_{i=1}^k
\mathbf{ w_i} = \mathbf w$). Moreover, let  $\bar R\in  \mathcal{R}(\mathcal H)$ be the region of $\Delta (\mathcal H)$
such that $c_{HM}(\mathcal H)=n(\mathcal H)- m(\bar R)$ ($n(\mathcal H)$ being the number of singular vertices of
$\mathcal H$ and $m(\bar R)$ being the number of singular vertices contained in $\bar R$). Note that $\bar R$ may be assumed to be the region of $\Delta (\mathcal H)$ obtained by ``fusing" the regions $\bar
R_1,$ $\bar R_2,$ $\dots,$ $\bar R_k$ ($\bar R_i$ being a suitable
region of $\Delta (\mathcal H_i)$ with $m(\bar R_i) \ne 0$
singular vertices, and $ \sum_{i=1}^k m(\bar R_i)= m(\bar
R)$).\footnote{Roughly speaking, we can say that $\bar R$ is the
``external" region of the planar realization of $\mathcal H$, and that $\bar R_i$ is the ``external" region of the planar
realization of $\mathcal H_i$, for each $i=1, \dots, k$.}
In fact, if this is not the case, it is easy to check that a new Heegaard diagram $\mathcal H^{\prime}$ of
$M$ with this property exists, with  $c_{HM}(\mathcal H^{\prime})
< c_{HM}(\mathcal H),$ against the hypothesis $c_{HM}(\mathcal
H)=c_{HM}(M)$.
Note that the diagram $\mathcal H_i$ ($i =1, \dots, k$) may fail to be reduced;
however, $c_{HM}(\mathcal H_i) \le n(\mathcal H_i)- m(\bar R_i)$ trivially holds by definition itself, for each $i=1, \dots, k.$

The thesis now directly follows:
$ \sum_{i=1}^k c_{HM}(\mathcal H_i) \le \sum_{i=1}^k n(\mathcal H_i)-  \sum_{i=1}^k m(\bar R_i) = n(\mathcal
H)- m(\bar R) = c_{HM}(\mathcal H).$
\qed

In order to prove our main result, we need a further lemma, which involves the singular 3-manifold $\widehat{M}$ obtained by capping off each component of $\partial M$ by a cone:

\begin{lemma}\label{singular_manifold} If $\mathcal H=(\Sigma_g, \mathbf {v}=\{v_1,v_2,\ldots, v_g\}, \mathbf {w}=\{w_1,w_2\ldots, w_s\})$ is a (reduced) Heegaard diagram of a compact $3$-manifold $M$ satisfying condition (*), then there exists a graph $\G \in G^{(0)}$ representing $\widehat{M}$ and regularly embedding into $\Sigma_g$, such that $c_{GM}(\G)=c_{HM}(\mathcal H)$.\end{lemma}

\dimo $\G$ is obtained by applying to $\mathcal H$ a construction described in \cite{[Cr]},
which is summarized briefly below:
\begin{itemize}
\item [-] Consider a planar realization of $\mathcal H$, where each curve of $\mathbf {v}$ (defining a 1-handle of $\Sigma_g$) is represented by two circles lying in the upper and lower half-plane respectively.  All circles in the upper half-plane may be assumed to have the same orientation; in the orientable (resp. non-orientable) case, all circles in the lower half-plane have the opposite orientation (resp. for each non-orientable 1-handle in $\Sigma_g$, the circle in the lower half-plane has the same orientation as its corresponding one in the upper half-plane).
    Let $v_0$ be the simple closed curve in $\Sigma_g$, represented by the $x$-axis in the plane. We can suppose that $v_0$ is disjoint from the curves of $\mathbf {v}$ and set  $\mathbf {v^\prime}= \mathbf {v}\cup \{v_0\}$.
Furthermore, by condition (*) and a suitable choice of $v_0$, we can always assume that the diagram  $(\Sigma_g, \mathbf {v^\prime}, \mathbf {w})$ is connected.
\item [-] Let $\mathbf {w^\prime}$ be the system obtained from $\mathbf {w}$ by ``doubling'' each of its curves.
For each $i=1,\ldots,s$, we denote by $\bar w_i$ the ``double'' of the curve $w_i$.
\item [-] Let $\G$ be the 4-coloured graph obtained by coloring alternatively 0 and 2 the arcs into which each curve of $\mathbf {v^\prime}$ is split by $\mathbf {w^\prime}$ and by coloring alternatively 1 and 3 the arcs into which each curve of $\mathbf {w^\prime}$ is split by $\mathbf {v^\prime}$. We assume to colour 0 the arcs lying between a curve of $\mathbf {w}$ and its ``double''. Furthermore we colour 1 (resp. 3) the arcs of $\mathbf {w^\prime}$ lying in the upper (resp. lower) half-plane of the planar realization.
\end{itemize}

Obviously,  $\mathcal H^\prime=(\Sigma_g,\mathbf {v^\prime},\mathbf {w^\prime})$ is a generalized Heegaard diagram of $M$, too.
Furthermore, $\G$ represents $\widehat M$ in virtue of \cite[Lemma 5]{[Cr]}).

By directly computing the Euler characteristics of the disjoint links of the vertices of $K(\G)$, it is easy to prove that all singular vertices of $K(\G)$ must be 0-labelled (see \cite[Lemma 4]{[Cr]}), i.e. $\G \in G^{(0)}.$
Hence, the equality $c_{GM}(\G)=c_{HM}(\mathcal H)$ follows by arguments of Section 5.
\qed

An example of the above construction is presented in Figures 6(a) and 6(b), starting from the planar realization of a genus two reduced
Heegaard diagram of $U_2\times I$, where all circles are assumed to have the same orientation.

 \bigskip

\centerline{\scalebox{0.8}{\includegraphics{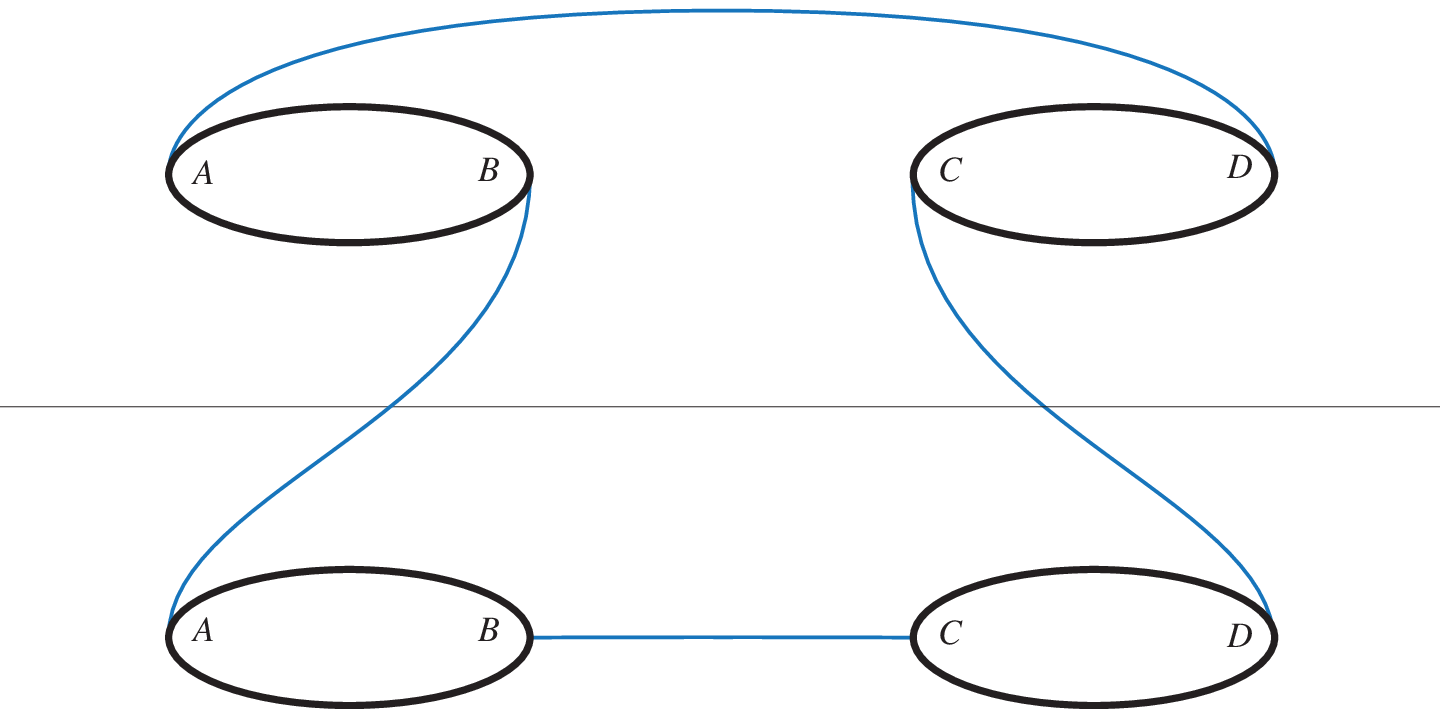}}}

\smallskip

\centerline{\footnotesize{Figure 6(a): $\mathcal H$, Heegaard diagram of $U_2 \times I$}}
\bigskip

\centerline{\scalebox{0.8}{\includegraphics{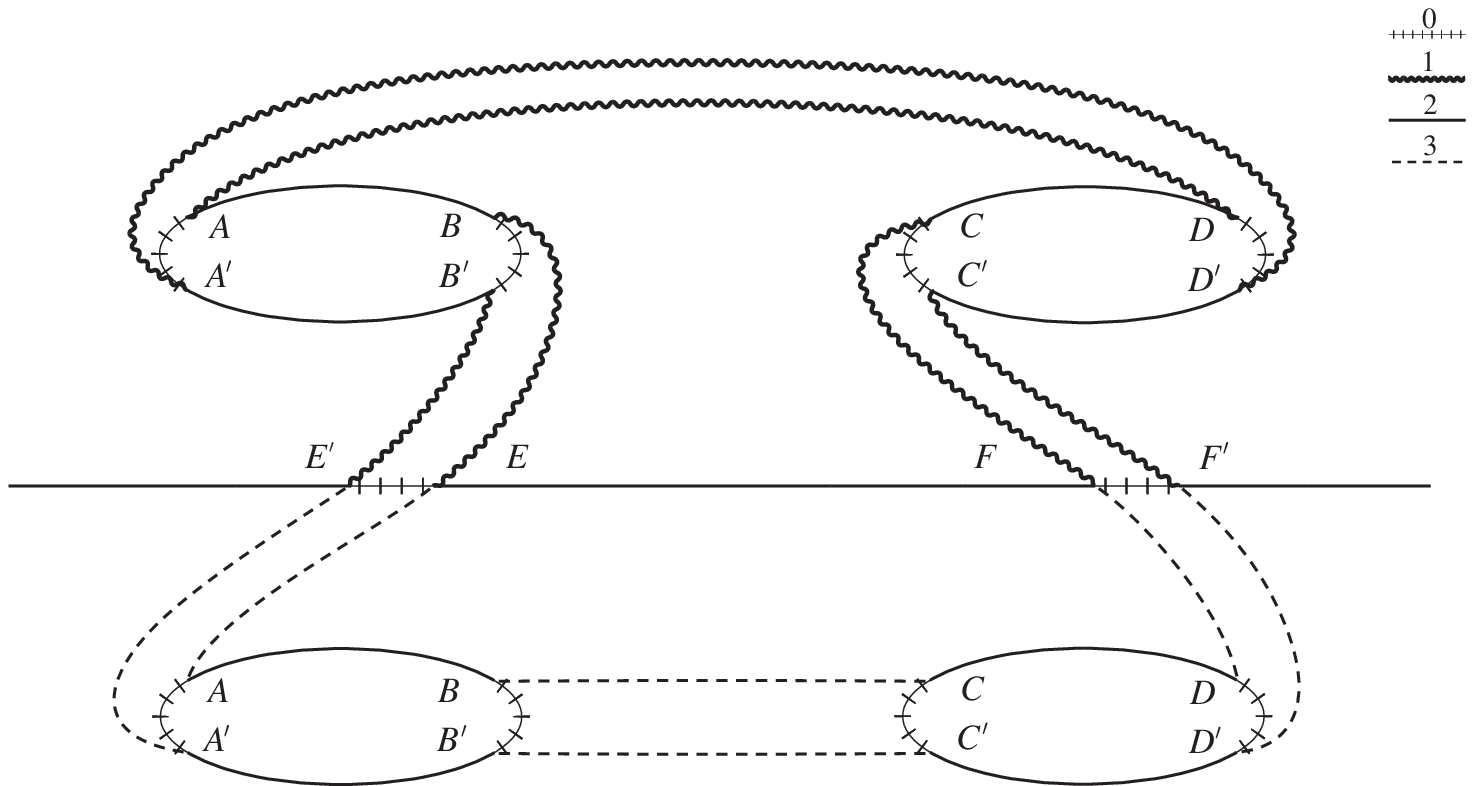}}}

\smallskip

\centerline{\footnotesize{Figure 6(b): $\G \in G^{(0)}$, representing $\widehat M$, with $M=U_2 \times I$}}

\bigskip
\bigskip

\begin{rem}\label{singular_vertices}\emph{Note that $K(\G_{\hat 0})$ can be obtained by cutting $\Sigma_g$ along the curves of $\mathbf {w}$. So, if $\mathbf {w}$ is reduced, each $\hat 0$-residue of $\G$ represents a surface of genus $>0$ and hence all $0$-labelled vertices of $K(\G)$ are singular.}\end{rem}

\bigskip
\bigskip

\textbf{Proof of the main result}
Let $M$ be a compact irreducible and boundary-irreducible 3-manifold and let $\mathcal H=(\Sigma_g,\mathbf {v}=$ $\{v_1,v_2,\ldots, v_g\},$ $\mathbf {w}=\{w_1,w_2\ldots, w_s\})$  be a reduced Heegaard diagram of $M$ such that $c_{HM}(\mathcal H)=c_{HM}(M)$.

\medskip
We will prove the statement by first assuming $\mathcal H$  connected  (case A), and then by taking into account a not connected  $\mathcal H$ (case B).

\medskip
\noindent {$\bullet$}  \  [{\bf case A}]  \ If $\mathcal H$  is assumed to be connected (and hence it satisfies condition (*)),
by applying the constructions of Lemma \ref{singular_manifold} and Lemma \ref{bis-tris} we obtain a gem $\bar\G$ of $M$.

Let us consider the subcomplex $K_{02}$ (resp. $K_{13}$) of $K(\Gamma)$ generated by the $0-$ and $2-$ (resp. $1-$ and $3-$) labelled vertices.
Note that the curves of the system $\mathbf {v^\prime}$ (resp. $\mathbf {w^\prime}$) are dual to the  edges of $K_{13}$ (resp. $K_{02}$).
It is obvious that, in order to reduce $\mathbf {v^\prime}$ (resp. $\mathbf {w^\prime}$) it is sufficient to remove the curve $v_0$ (resp. the curve $\bar w_i$,  i.e. the ``double'' of the curve $w_i$ of $\mathbf {w}$, for each $i=1,\ldots,s$).

\medskip

\centerline{\scalebox{0.65}{\includegraphics{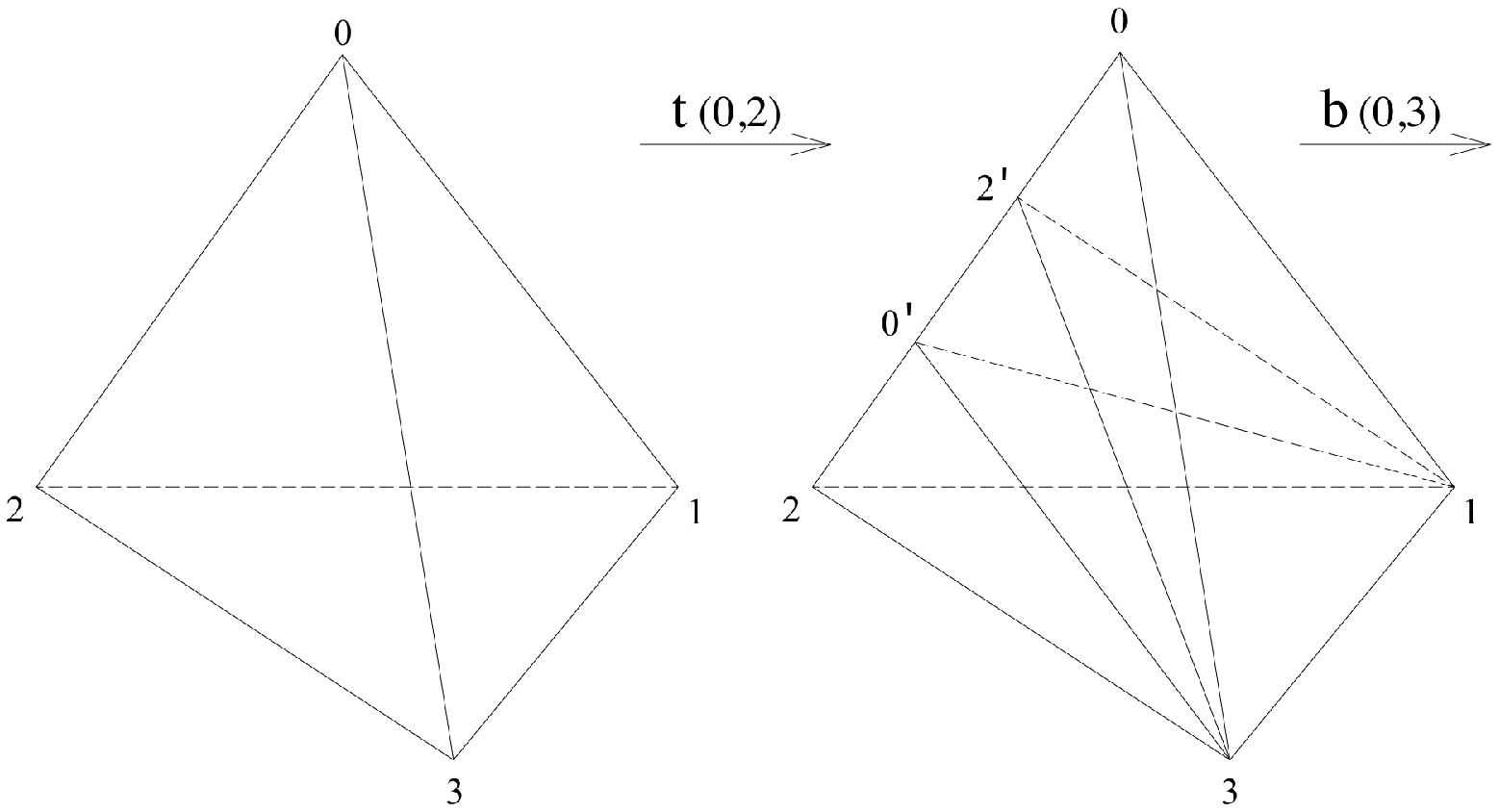}}}
\smallskip

\bigskip
\medskip

\centerline{\scalebox{0.75}{\includegraphics{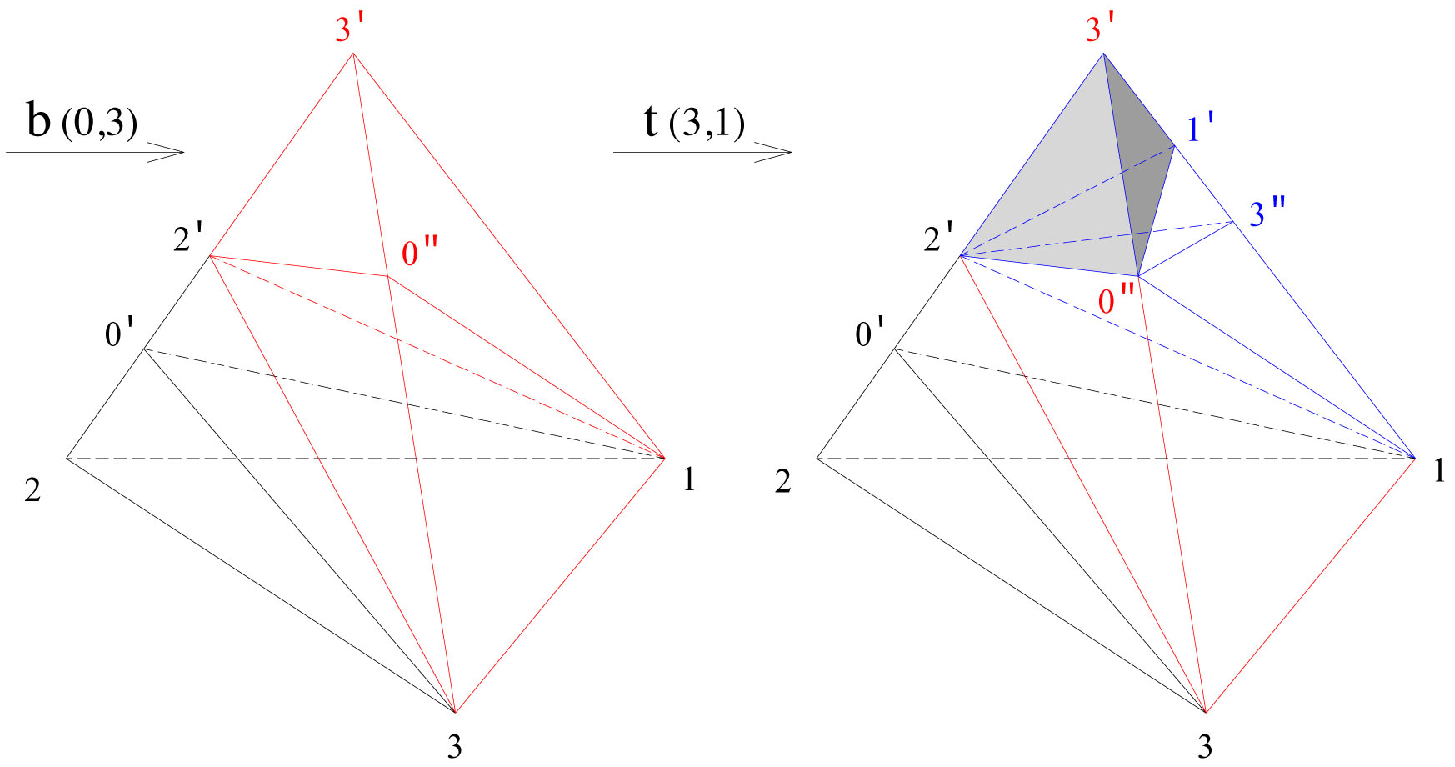}}}

\smallskip

\centerline{\footnotesize{Figure 7}}

\bigskip

We turn now our attention to the pseudocomplex $\bar K$.
Figure 7 describes locally, in a tetrahedron near a 0-labelled singular vertex of $K(\Gamma)$, the effects of the bisection and trisections performed. In the last picture of the sequence, we  highlighted the small tetrahedron which has to be removed in the last step of the construction, as part of a neighborhood of the singular vertex.
Moreover, we put ``primes"  beside the colour for the new vertices arising from the subdivisions.

In the following we will call an edge of $K(\bar\Gamma)$ \textit{of type $i-j\ $} (with $i,j\in\Delta_3\cup\{0^\prime,0^{\prime\prime},1^\prime,2^\prime,3^\prime,3^{\prime\prime}\}$) referring to the labels of its endpoints as in Figure 7.

If we consider the generalized Heegaard diagram $\mathcal H_\alpha=(\Sigma^{(\alpha)},\mathbf {v_\alpha},\mathbf {w_\alpha})$ associated to $\bar\Gamma$ with respect to the choice $\alpha=1$, by Lemma \ref{bis-tris}, we have $\Sigma^{(\alpha)}=\Sigma_g$.

The system $\mathbf {v_\alpha}$ is formed by the curves of $\mathbf {v^\prime}$ and by the curves dual to the edges of type $1-3^{\prime\prime}$ and $1^\prime-3^{\prime\prime}$.
The system $\mathbf {w_\alpha}$ is formed by the curves dual to the edges of type $0^\prime-2$ and $0^\prime-2^\prime$, which correspond in pairs to the curves of $\mathbf {w^\prime}$.

Let us call $\bar{\mathcal D}$ the set formed by $v_0$ and by the curves of $\mathbf {v_\alpha}$ which are dual to the edges of type $1^\prime-3^{\prime\prime}$ and $1-3^{\prime\prime}$.
$\bar{\mathcal D}$ corresponds to a maximal tree of $\bar K_{13}$, since $v_0$ is dual to a maximal tree of  $K_{13}$ and no edges of type $1^\prime-3^{\prime\prime}$ have a common endpoint in $\bar K$.

Let $G_{02}$ be the graph obtained by $\bar K_{02}$ by shrinking to a point each connected component of $\partial\bar K$, and let $\bar G$ be the subgraph of $G_{02}$ formed by all edges of type $0^\prime-2^\prime$ (modulo the above identifications) and by those edges of type $0^\prime-2$ belonging to the subdivision of an edge of $K$ dual to a curve $\bar w_i$ ($i=1,\ldots,s$)\footnote{Note that, by Remark \ref{singular_vertices}, all edges of $K_{02}$ have been subdivided.}.

We point out that:
\begin{itemize}
\item[-] by Remark \ref{singular_vertices}, each vertex of $G_{02}$ belongs to $\bar G$;
\item[-] each connected component of $\bar G$ is a tree; in fact, the edges of $K(\G)$ dual to the curves $\bar w_i$ form a maximal tree of $K_{02}$, no edges of type $0^\prime-2^\prime$ have the $0$-labelled endpoint in common and no edges of type $0^\prime-2$ belonging to $\bar G$ have a common endpoint;
\item[-] each connected component of $\bar G$ intersects only one boundary component of $\bar K$ (otherwise the edges of $K(\G)$ dual to the curves $\bar w_i$ wouldn't form a tree of $K_{02}$).
\end{itemize}

If we denote by $\bar{\mathcal D}^\prime$ the set of curves of $\mathbf {w_\alpha}$ dual to the edges of $\bar G$, as a consequence of the above considerations and of Proposition \ref{reduced Heegaard}, we have that $\mathcal H_\alpha (\bar{\mathcal D},\bar{\mathcal D}^\prime)$ is a reduced Heegaard diagram of $M$ associated to $\mathcal H_\alpha$; furthermore it is obvious that $\mathcal H_\alpha (\bar{\mathcal D},\bar{\mathcal D}^\prime)=\mathcal H$.

\medskip

Hence, $c_{GM}(\bar\G)=c_{HM}(\mathcal H_\alpha)=c_{HM}(\mathcal H)$, as claimed.

\bigskip

\noindent {$\bullet$}  \  [{\bf case B}]  \  Let us now assume that the reduced Heegaard diagram $\mathcal H=(\Sigma_g,\mathbf {v}=\{v_1,v_2,\ldots, v_g\},$ $\mathbf {w}=\{w_1,w_2\ldots, w_s\}),$ with the property $c_{HM}(\mathcal H)=c_{HM}(M)$, is not connected.

By definition, this means that at least a connected component of $\Sigma_g\setminus (\mathbf {v}\cup \mathbf {w})$ is not a  2-disk.
As a consequence, a simple closed essential curve $l$ exists in $\Sigma_g$, not intersecting $\mathbf {v}\cup \mathbf {w}.$  Obviously, $l$ bounds a disk $D$ in the handlebody $\mathbb X_g$ of the Heegaard splitting $(\Sigma_g, \mathbb X_g, \mathbb Y_g)$ associated to $\mathcal H$; on the other hand, in the compression body $\mathbb Y_g$, $l$ bounds either a disk $D^\prime$ or an annulus $A$ with $\partial A$ consisting of the disjoint union of $l$  ($l \in \partial_{+} \mathbb Y_g$)  and $l^\prime$ ($l^\prime \in \partial_{-} \mathbb Y_g=\partial M$).

In the first case, i.e. if a disk $D^\prime \subset \mathbb Y_g$ exists with $\partial D^\prime =l$, then the union $D \cup D^\prime$ yields a 2-sphere in $M,$ which - because of the assumption of irreducibility - splits $M$ into a (necessarily trivial) connected sum. This means  that $l$ splits $\mathcal H$ into two Heegaard diagrams $\mathcal H^{\prime}$ and $\mathcal H^{\prime \prime}$, where $\mathcal H^{\prime}$ represents $M,$ too, and $\mathcal H^{\prime \prime}$ represents $\mathbb S^3$  (or viceversa).

In the second case, i.e. if an annulus $A \subset \mathbb Y_g$ exists with $\partial A$ consisting of $l \subset \partial_{+} \mathbb Y_g$
and $l^\prime \subset \partial_{-} \mathbb Y_g=\partial M$,  then $\bar D = D \cup A$ is a proper disk in $M,$  with $\partial \bar D = l^\prime \subset \partial M,$ which - because of the assumption of boundary-irreducibility -  splits $M$ into a (necessarily trivial) boundary connected sum.
This means  that $l$ splits $\mathcal H$ into two Heegaard diagrams $\mathcal H^{\prime}$ and $\mathcal H^{\prime \prime}$, where $\mathcal H^{\prime}$ represents $M,$ too, and $\mathcal H^{\prime \prime}$ represents $\mathbb D^3$  (or viceversa).

By iterating the process for any possible subsequent splitting of $\mathcal H,$ we obtain $k \ge 2$ connected Heegaard diagrams  $\mathcal H_1,$  $\mathcal H_2,$ $\dots,$ $\mathcal H_k$ ($k \ge 2$) so that  $\mathcal H_1$ represents the compact 3-manifold $M$, too,  while for each $i=2, \dots, k,$  $\mathcal H_i$ represents either the 3-sphere or the 3-disk.

Now, Lemma \ref{lemma2} ensures
$\sum_{i=1}^k c_{HM}(\mathcal H_i) \le c_{HM}(\mathcal H),$ and so $c_{HM}(\mathcal H_1) \le c_{HM}(\mathcal H)$ directly follows.

Since $c_{HM}(\mathcal H)=c_{HM}(M)$ is assumed,  $c_{HM}(\mathcal H_1) = c_{HM}(\mathcal H)= c_{HM}(M)$ necessarily holds, too.

Hence, the general statement follows by applying the first part of the proof (case A) to any connected reduced Heegaard diagram  of $M$ associated to $\mathcal H_1$.
\qed

\medskip

\begin{rem}\label{condizione_stella}\emph{Note that, since Lemma \ref{singular_manifold} holds for any Heegaard diagram $\mathcal H$ satisfying condition (*), the proof of the main theorem (i.e. the coincidence between $GM$-complexity and $HM$-complexity) works as well for all compact 3-manifolds admitting such a Heegaard diagram. More precisely:} if $M$ is a compact 3-manifold admitting a Heegaard diagram $\mathcal H$ satisfying condition (*) such that $c_{HM}(\mathcal H)=c_{HM}(M)$, then a gem $\bar \G$ of $M$ may be directly obtained, with the property $c_{GM}(\bar\G)=c_{HM}(\mathcal H)= c_{HM}(M)$. \end{rem}

\bigskip

Moreover, as a consequence of the main result, we can establish:

\begin{prop} \label{uguaglianza complessità varieta' singolari}

For each compact irreducible and boundary-irreducible 3-manifold $M$, $\hat c_{GM}(M)= c_{GM}(M)$.\end{prop}

\dimo\ Let $\G^\prime\in G^{(0)}$ represent $\widehat{M}$, with  $c_{GM}(\G^\prime)=\hat c_{GM}(M)$. Then, there is a reduced Heegaard diagram $\mathcal H_\alpha$ of $M$, associated to $\G^\prime$, such that $c_{GM}(\G^\prime)=c_{HM}(\mathcal H_\alpha)$; moreover, by the hypothesis on $M$, the proof of the main result (case B) allows to assume $\mathcal H_\alpha$ connected, without loss of generality.
Hence, by applying to $\mathcal H_\alpha$ the procedure described in the proof of the main theorem (case A), a 4-coloured graph $\G,$ regular with respect to colour $3$ and representing $M$ is obtained. Since $\G$ has the same GM-complexity as $\mathcal H_\alpha$, $c_{GM}(M)\leq\hat c_{GM}(M)$ follows.

Conversely, let $\G$ be a 4-colored graph with boundary representing $M$ and let $\mathcal H_\alpha$ be a reduced Heegaard diagram associated to $\G$ such that $c_{GM}(\G)=c_{HM}(\mathcal H_\alpha)$.
Again by the arguments of case B in the proof of the main theorem, we can suppose $\mathcal H_\alpha$ to be connected.
Therefore, by applying to $\mathcal H_\alpha$ the construction of Lemma \ref{singular_manifold}, we get a 4-coloured graph $\G^\prime$, representing $\widehat{M}$, such that $c_{GM}(\G^\prime)=c_{HM}(\mathcal H_\alpha)$. Hence $\hat c_{GM}(M)\leq c_{GM}(M)$.
\qed

\bigskip

In the following section we will apply  the notion of GM-complexity (in particular, via graphs representing singular 3-manifolds) in order to yield an estimation of Matveev's complexity for an infinite family of bounded 3-manifolds, including all torus knot complements.

\bigskip
\bigskip

\section{\hskip -0.7cm . Estimation of Matveev's complexity for torus knot complements}

In \cite[Figure 1d]{[Gr]}, Grasselli obtained a 4-coloured graph $ \bar \Lambda$ representing a singular manifold which - in virtue of a previous work  by Montesinos (see  \cite[Example 2]{[Mo]}) - turns out to be associated to the complement of the trefoil knot $\frak t(3,2)$ in $\mathbb S^3$. By a suitable combinatorial move not affecting the represented manifold (i.e., by the inverse of a $\rho_2$-pair switching, which may be performed via the insertion of a 2-dipole, followed  by the elimination of a (proper) 1-dipole: see \cite{[BCG]}
for details on the admissible combinatorial moves on edge-coloured graphs representing PL-manifolds and \cite{[G$_1$]}, \cite{[BC]}  for their extension to the case of singular manifolds), the 4-coloured graph $ \bar \Lambda^{\prime}$  of Figure 8(a) is obtained.

\bigskip

The definitions and results of the previous sections easily yield a direct proof that the trefoil knot complement has Matveev's complexity zero:\footnote{The fact is actually already known, as Sergei Matveev pointed out in a personal communication, referring to an unpublished work by his former student Nikolaev, which yields the complete classification of all compact irreducible and boundary irreducible 3-manifolds with complexity zero.}

\begin{prop} \ $ c( \mathbb S^3 - \frak t(3,2)) =0.$ \end{prop}

\dimo \ In virtue of Proposition \ref{complessità_singolare}, Matveev's complexity may be estimated from any 4-coloured graph representing the associated singular manifold, and having only one singular colour. A trivial check ensures that the graph $ \bar \Lambda^{\prime}$ actually belongs to the class $G^{(0)}$: in fact, all $\hat c$-residues represent 2-spheres, except when $c=0$ (the Euler characteristic computation easily proves $ \bar \Lambda^{\prime}_{\hat 0}$ to represent a torus). Moreover, since $g_{\hat c}=1 \ \forall c\in \Delta_3$, one only bicoloured cycle of each complementary pair of colours has to be deleted, in order to obtain a reduced Heegaard diagram of the trefoil knot complement.  Hence,  $c_{GM}(\bar \Lambda^{\prime}) = 5- 5=0$ directly follows: see Figure 8(b) for an associated Heegaard diagram, where it is easy to detect a region containing all the intersection points of the two systems of curves.
\ \ \ \ \ \qed

\medskip

\centerline{\scalebox{0.85}{\includegraphics{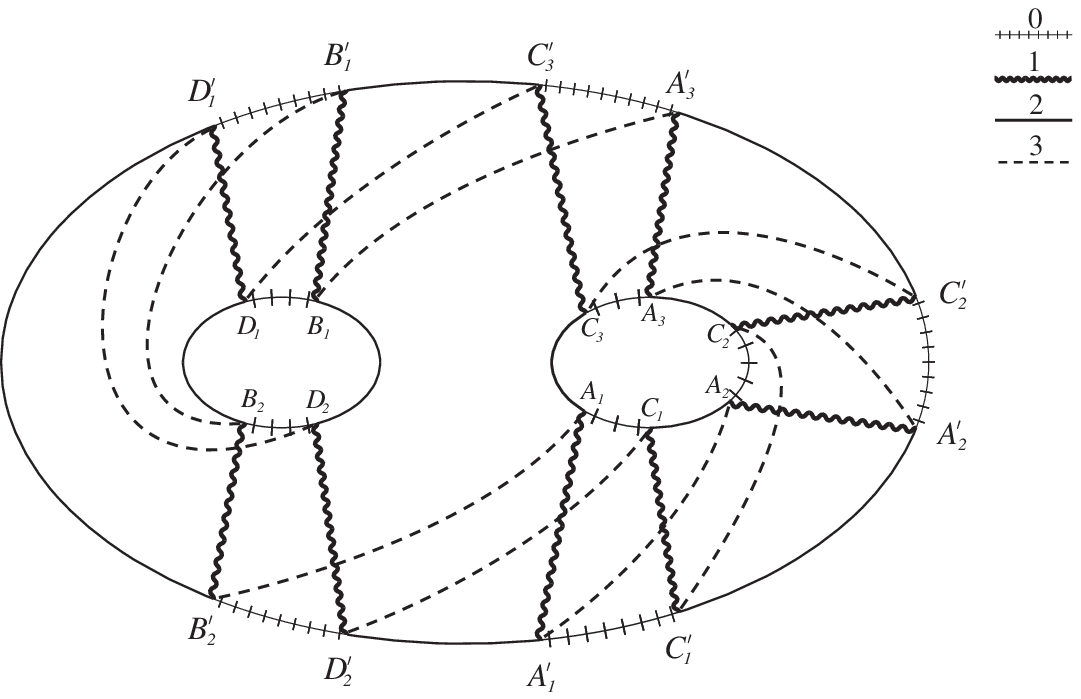}}}

\smallskip

\centerline{Figure 8(a):  $ \bar \Lambda^{\prime}$, representing the trefoil knot complement}

\bigskip

\medskip

\centerline{\scalebox{0.6}{\includegraphics{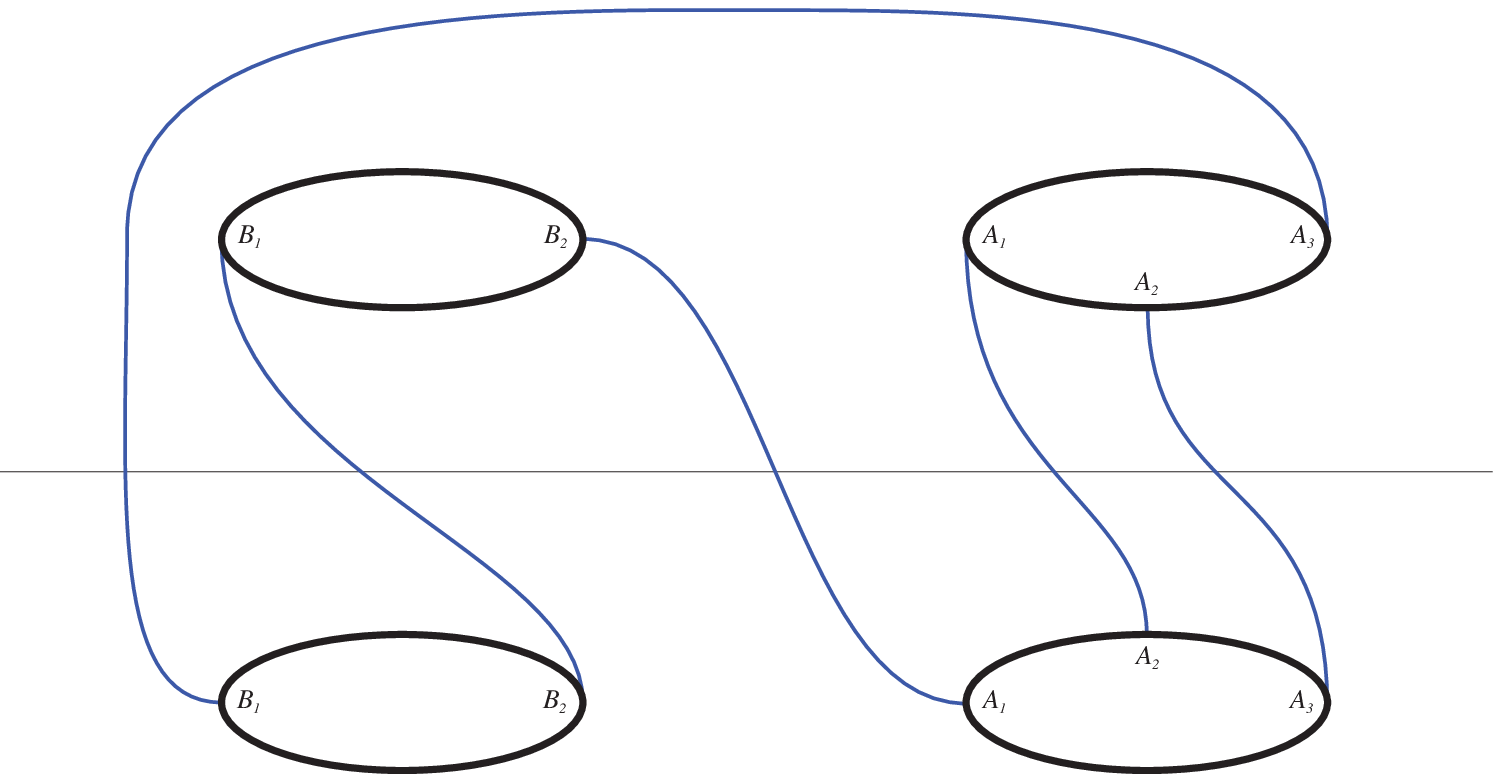}}}

\smallskip

\centerline{\footnotesize{Figure 8(b): a reduced Heegaard diagram of the trefoil knot complement}}

\bigskip
\medskip

By a quite natural generalization of the 4-coloured graph $ \bar \Lambda^{\prime}\in G^{(0)}$ representing the trefoil knot complement, let us now consider the 4-coloured graph $\Lambda((p,h),(q,k))$ depicted in Figure 9(a): it depends on two pairs $(p,h)$ and $(q,k)$ of coprime integers (where $p \ge h \ge 1$, $q \ge k \ge 1$ and $p \ge q \ge 2$ may be assumed).

\bigskip

\centerline{\scalebox{0.85}{\includegraphics{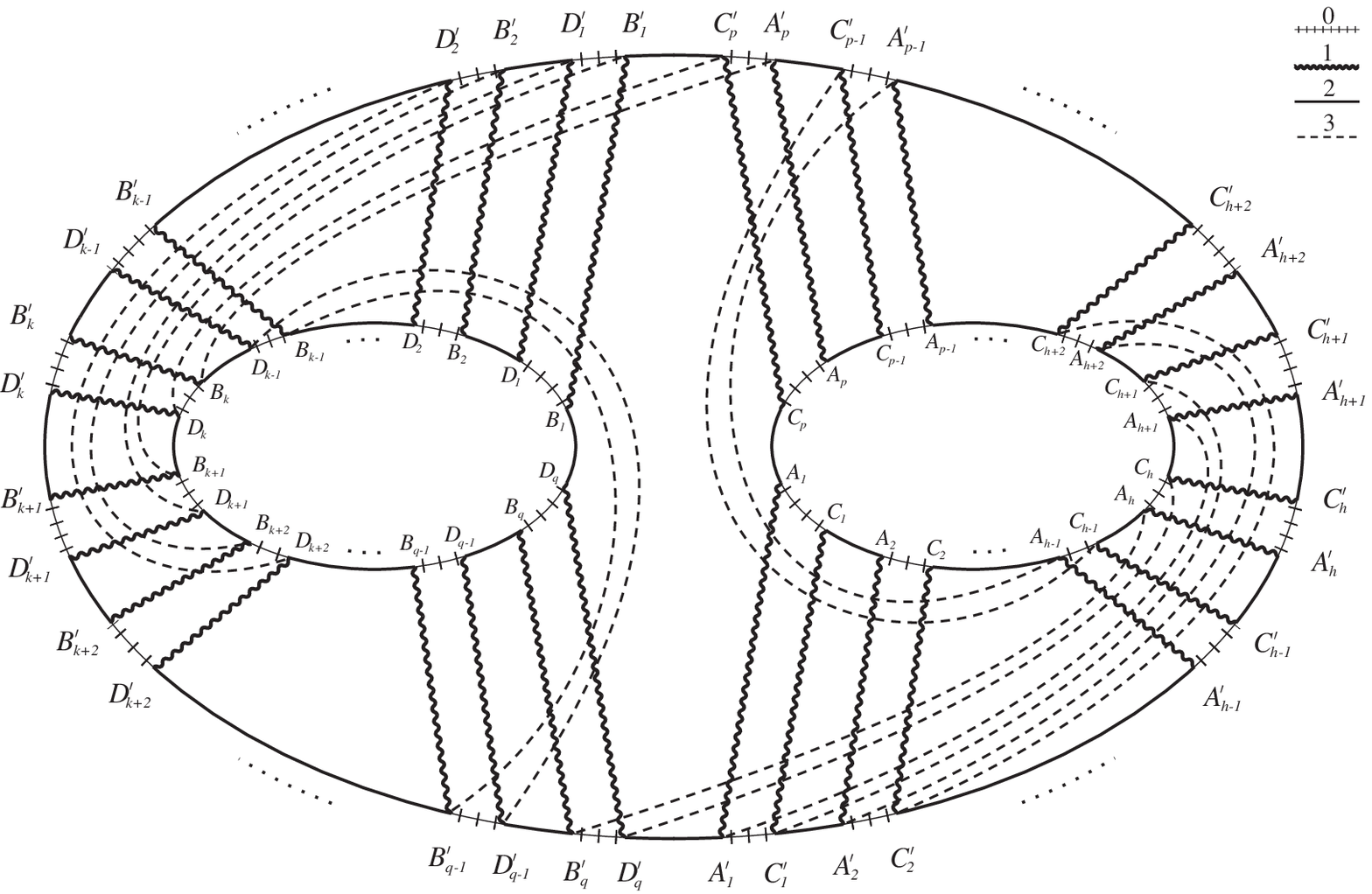}}}

\smallskip

\centerline{\footnotesize{Figure 9(a): $\Lambda((p,h),(q,k))$}}

\bigskip

For any 4-tuple $(p,h,q,k)$ of positive integers, 3-coloured edges may be described via the following pairs of adjacent vertices:
\begin{itemize}
\item{} $\forall i = 1, \dots, p-1:$ $(A_i^\prime, A_{i+h})$ and $(C_i^\prime, C_{i+h})$, where indices are assumed $\mod p$ within the set $\{1,2, \dots, p\};$
\item{} $\forall j = 1, \dots, q-1:$ $(B_j^\prime, B_{j+k})$ and $(D_j^\prime, D_{j+k})$, where indices are assumed $\mod q$ within the set $\{1,2, \dots, q\};$
\item{} $(A_p^\prime, B_k)$ and $(C_p^\prime, D_k)$; \ $(B_q^\prime, A_h)$ and $(D_q^\prime, C_h)$.
\end{itemize}

\medskip

\begin{prop} \label{identificazione Seifert} For any 4-tuple $(p,h,q,k)$ of positive integers so that $GCD(p,h)=GCD(q,k)=1$, the 4-coloured graph $\Lambda((p,h),(q,k))$ belongs to the class $G^{(0)}$. Moreover, it represents the Seifert 3-manifold $(\mathbb D^2; (p,\alpha),(q,\beta))$, with base $\mathbb D^2$ and two exceptional fibers, where the Seifert type $(p, \alpha)$ (resp. $(q,\beta)$)  of the first (resp. second) fiber is uniquely determined by $\alpha h \equiv 1 \mod p$ (resp. $\beta k \equiv 1 \mod q$).
\end{prop}

\dimo \ First of all, note that the 4-coloured graph $\Lambda((p,h),(q,k))$ (where $GCD(p,h)=$ $GCD(q,k)=1$ is assumed to hold) has the following combinatorial structure: $ g_{01}= p+q;$ \ $g_{02}= 3;$ \ $g_{03}= p+q;$ \ $g_{12}= p+q-1;$ \ $g_{13}= 2;$  \ $g_{23}= p+q-1.$

So the regularity property, together with the connectedness of each 3-coloured subgraph, allows to topologically recognize $K(\Gamma_{\hat c})$ \ ($\forall c \in \Delta_3$) by means of an Euler characteristic computation:
$$ \chi (K(\Gamma_{\hat 0})) = 0  \ \ \ \Rightarrow \ \ \ K(\Gamma_{\hat 0}) \ \text {is \ a \ torus}$$
$$ \chi (K(\Gamma_{\hat 1})) = \chi (K(\Gamma_{\hat 2})) = \chi (K(\Gamma_{\hat 3}))  = 2  \ \ \ \Rightarrow \ \ \ K(\Gamma_{\hat c})  \ \ \text {is \ a \ 2-sphere} \ \ \forall c \in \{1,2,3\}$$
This proves that $\Lambda((p,h),(q,k))$ belongs to the class $G^{(0)}$ and that the represented compact 3-manifold has (connected) toric boundary.

Finally, in order to complete the proof, it is sufficient to take into account the reduced Heegaard diagram - $\mathcal H ((p,h),(q,k)),$ say - associated to $\Lambda((p,h),(q,k))$ with respect to colour $2$, where the subset $\mathcal D$ (resp. $\mathcal D^{\prime}$) considered in Section 5 consists of the length $2(p+q)$ $\{0,2\}$-coloured cycle (resp. of one arbitrarily chosen $\{1,3\}$-coloured cycle): see Figure 9(b). Since $\mathcal H((p,h),(q,k))$  is a simple closed curve on the genus two handlebody consisting of exactly two connections of type $(p,h)$ and $(q,k)$  (see \cite{[Z]}), it turns out to be isotopic to the so-called ``standard Heegaard diagram $HD_0$ of type $(S^p T^{-q}; h, k)$" described in \cite[Definition 4.3]{[BRZ]}.

\medskip

\centerline{\scalebox{0.7}{\includegraphics{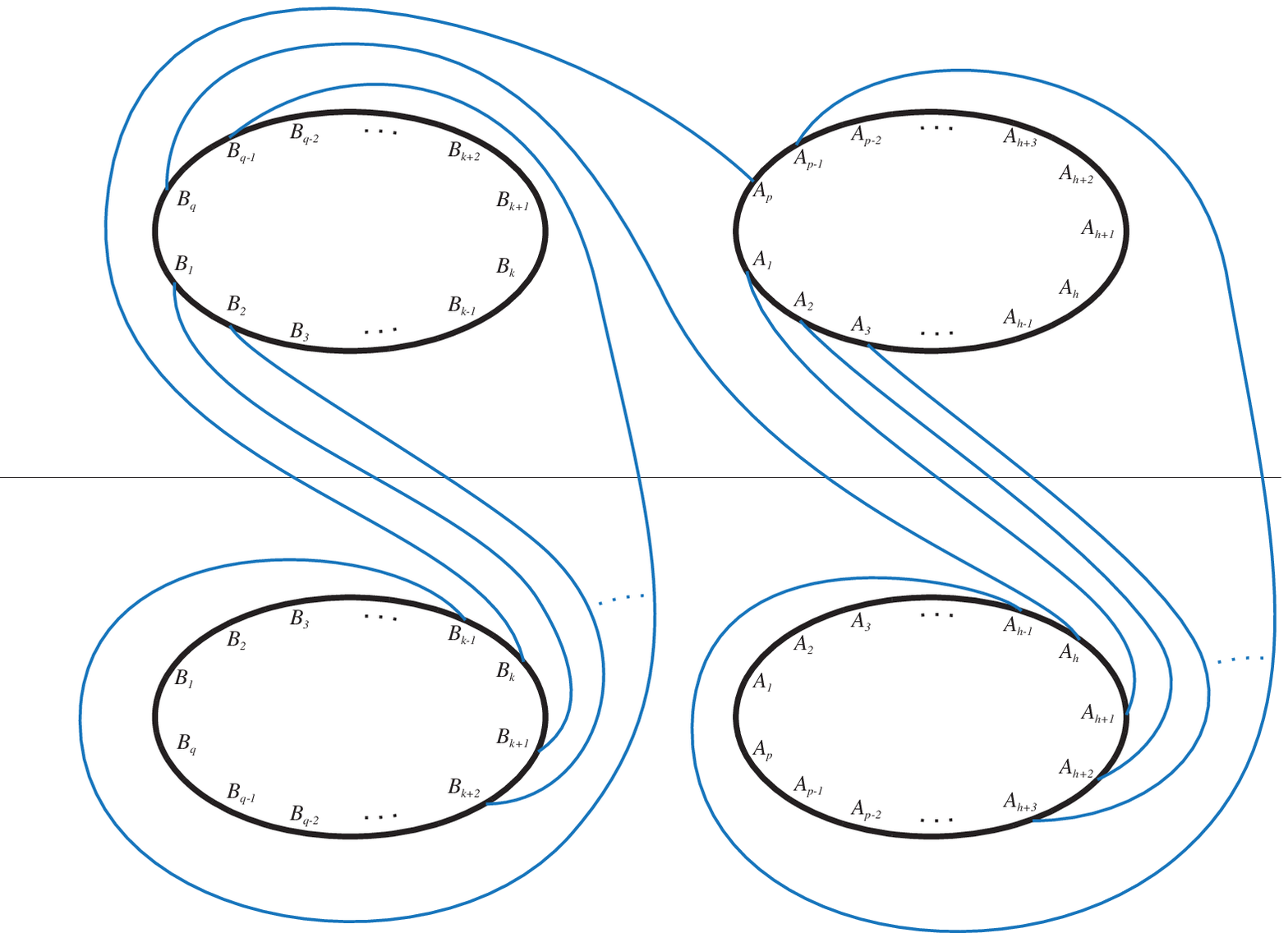}}}

\smallskip

\centerline{\footnotesize{Figure 9(b): $\mathcal H ((p,h),(q,k))$}}

\bigskip

The statement now follows from \cite[Proposition 4.4]{[BRZ]}: the standard Heegaard diagram $HD_0$ of type $(S^p T^{-q}; h, k)$ is proved to represent the Seifert 3-manifold $(\mathbb D^2; (p,\alpha),(q,\beta))$, with $\alpha h \equiv 1 \mod p$ and $\beta k \equiv 1 \mod q$, and hence both $\mathcal H ((p,h),(q,k))$  and $\Lambda((p,h),(q,k))$ do, too.
 \ \ \ \ \
 \qed

\vskip-0.6truecm

\begin{rem}\label{remark_raddoppio}\emph{It is worthwhile noting that $\Lambda((p,h),(q,k))$ may be directly obtained from $\mathcal H ((p,h),(q,k))$ by applying the construction described in Lemma \ref{singular_manifold}.} \end{rem}

\smallskip

\begin{cor} For each pair $(p,q)$ of coprime integers ($p > q \ge 2$),  $\Lambda((p,q),(q,p))$ represents the complement of the torus knot $\frak t(p,q)$.
\end{cor}

\dimo  \ By \cite[Corollary 4.6]{[BRZ]}, the complement of the torus knot $\frak t(p,q)$ (which is the Seifert manifold $(\mathbb D^2; (p,\alpha),(q,\beta))$, with  $ \alpha \cdot q \equiv \pm 1 \mod p$ and $ \beta \cdot p \equiv \pm 1 \mod q$, by  \cite[Proposition 4.2 (d)]{[BRZ]}) admits the standard Heegaard diagram $HD_0$ of type $(S^p T^{-q}; q, p)$. Hence, the thesis is a direct consequence of the proof of the previous Proposition: $\mathbb S^3 - \frak t(p,q)$ is also represented by the isotopic Heegaard diagram $\mathcal H ((p,q),(q,p)),$
as well as  by the graph  $\Lambda((p,q),(q,p))$  (i.e. $\Lambda((p,q),(q,p^{\prime}))$, where $1 \le p^{\prime} \le q-1$ is (uniquely) determined  by  $p^{\prime}\equiv p  \mod q$).
\ \ \ \ \
\qed

\vskip-0.6truecm

\begin{rem}\label{remark_curvarossa}\emph{As kindly pointed out by Jos\'e Maria Montesinos in a personal communication, a direct proof exists of the fact that the 4-coloured graph $\Lambda((p,q),(q,p))$  represents the complement in $\mathbb S^3$ of a knot (and hence of the appropriate torus knot, because of the fundamental group computation): it is sufficient to note that a simple closed curve (whose existence relies on the Bezout equation for the pair $(p,q)$ of coprime integers) may be added to the associated Heegaard diagram $\mathcal H ((p,q),(q,p))$ so that the represented manifold (i.e. the manifold obtained from $K(\Lambda((p,h),(q,k)))$ by a 2-handle addition) has null fundamental group.  See Figure 10 for an example, in the case $p=5$, $q=3$ (where the red dotted curve identifies the ``additional" 2-handle).}  \end{rem}

\centerline{\scalebox{0.65}{\includegraphics{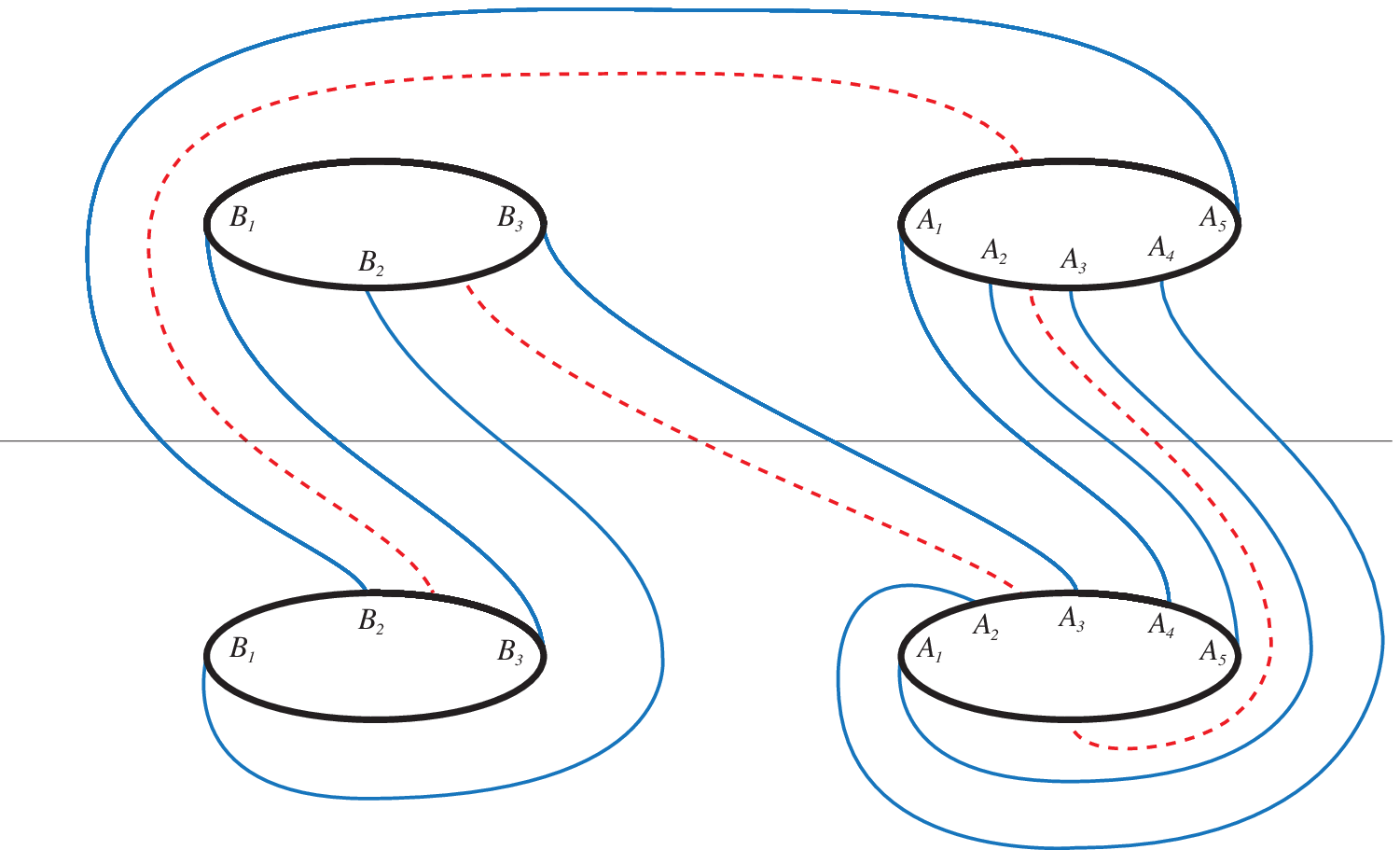}}}

\smallskip

\centerline{\footnotesize{Figure 10: $\mathcal H ((5,3),(3,2))$}}

\bigskip

As a consequence, an estimation for the Matveev's complexity of the represented manifolds is obtained via the notion of GM-complexity; in the case of torus knot complements it significantly improves, for this particular class of knots, the general estimation for Matveev's complexity of any link complement given in \cite[Prop. 2.11]{[Ma]}.\footnote{While writing the present paper, we learnt about a preprint by Fominykh and Wiest, which also yields an upper estimation for Matveev's complexity of torus knot complements, via their representation as Seifert manifolds: see \cite{[FW]}. In some cases, Fominykh and Wiest's estimation further improves Matveev's one.}

\begin{prop}  \label{complessità grafi parametrici}

Let \ $ \delta_{\alpha}= \begin{cases} 1 & \alpha \equiv \pm 1 \mod p \\
                                              0 & otherwise  \end{cases}$
                                                                                       \     and  \ $ \delta_{\beta}= \begin{cases} 1 & \beta \equiv \pm 1 \mod q \\
                                              0 & otherwise  \end{cases}$

\noindent Then:
$$ c((\mathbb D^2; (p,\alpha),(q,\beta)) \le
 \max \{p-4 + \delta_{\alpha}, 0\} + \max \{q-4 + \delta_{\beta}, 0\}.$$
In particular:
\begin{itemize}
\item[(a)]  \  For each pair $(p,q)$ of coprime integers ($p >  q > 3$), with $p-q \ne 1$,
$$ c( \mathbb S^3 - \frak t(p,q)) \le  p+q -8; $$
\item[(b)] \ for each $p \ge 4,$ \ $ c( \mathbb S^3 - \frak t(p,p-1)) \le 2p-7;$
\item[(c)] \ for each $p\ge 5$,  \ $ c( \mathbb S^3 - \frak t(p,2)) \le p-4;$
\item[(d)] \ for each $p \ge 5$,  \ $ c( \mathbb S^3 - \frak t(p,3)) \le p-4.$
 \end{itemize}
 \end{prop}

\dimo \  \
By arguments of Section 5, Matveev's complexity of $M= (\mathbb D^2; (p,\alpha),(q,\beta))$ (where $GCD(p,\alpha)=GCD(q,\beta)=1$ necessarily holds) may be estimated from any graph belonging to the class $G^{(0)}$ and representing the associated singular 3-manifold $\widehat M$ (and $M$ itself, too).
The general statement follows by considering the above described 4-coloured graph $\Lambda((p,h),(q,k))$, with $h,k$ uniquely determined by $\alpha h \equiv 1 \mod p$ and $\beta k \equiv 1 \mod q$, and by computing its $GM$-complexity by means of the associated reduced Heegaard diagram $\mathcal H((p,h),(q,k))$: there are exactly $p+q$ intersection points between the two systems of curves (i.e.: $\{A_1, A_2, \dots, A_p\}\cup \{B_1, B_2, \dots, B_q\}$), while a ``maximal" region of the embedding contains the subset  $\{A_1, A_p,$ $A_h,$ $A_{h+1}\}$ $ \cup$ $ \{B_1, B_q, B_k, B_{k+1}\}$,  whose cardinality depends on the value of $p$ and $q$  (since all indices are assumed to be $\mod p$ (resp. $\mod q$) within the sets $\{1,2, \dots, p\}$ (resp. $\{1, 2, \dots, q\}$)), as well as on the possible coincidence among elements.
In particular, coincidences occur in $\{A_1, A_p,  A_h, A_{h+1}\}$ (resp. in $\{B_1, B_q, B_k, B_{k+1}\}$) if and only if $ h \equiv \pm 1 \mod p$  (resp. $k \equiv \pm 1 \mod q$), i.e. if and only if   $\alpha \equiv \pm 1 \mod p$ (resp.  $\beta \equiv \pm 1 \mod q$).

Statements (a), (b), (c), (d) are easy consequences of the general relation.
\ \ \ \ \
\qed

\smallskip

The following Corollary collects some results which directly follow from Proposition \ref{complessità grafi parametrici}. However, they appear to be of particular interest, since the exact value of the Matveev's complexity is obtained, for all involved 3-manifolds: the complements of torus knots $\frak t(4,3)$, $\frak t(5,2)$ and $\frak t(5,3)$ are proved to have Matveev's complexity less or equal to one \footnote{These estimations actually turn out to be the exact values of complexity, since these complements of torus knots do not belong to the family of complexity zero compact 3-manifolds, according to Matveev's communication about Nikolaev's result.}, while the Seifert manifold $(\mathbb D^2; (2,1),(2,1))= U_2 \tilde \times I$ (i.e. the orientable I-bundle over the Klein Bottle) and the  Seifert manifold $(\mathbb D^2; (3,1),(3,1))$ are proved to have Matveev's complexity zero.\footnote{Note that the last two facts are already known, too, via Nikolaev's result.}

\begin{cor}  \label{complessità torus knot specifici}
\par \noindent
\begin{itemize}
\item[(a)] \ $ c(\mathbb S^3 - \frak t(4,3)) \le 1;$
\item[(b)] \ $ c(\mathbb S^3 - \frak t(5,2)) \le 1;$
\item[(c)] \ $ c(\mathbb S^3 - \frak t(5,3)) \le 1;$
\item[(d)] \ $ c(U_2 \tilde \times I)=0;$
\item[(e)] \ $ c((\mathbb D^2; (3,1),(3,1))=0.$
    \end{itemize}
 \end{cor}

\vskip-0.8truecm
\ \qed

\begin{rem}\label{caso particolare}\emph{Case (c) of the previous Corollary may be directly checked by the above Figure 10, where it is very easy to see a region of the embedding containing all intersection points, except $A_2$.} \end{rem}

\medskip

{\it Acknowledgement.} Work performed under the auspices of the
G.N.S.A.G.A. of I.N.d.A.M. (Italy) and financially supported by
M.I.U.R. of Italy, University of Modena and Reggio Emilia, funds for selected research topics.

\medskip

\footnotesize{
}

\end{document}